\documentclass[letterpaper, 10 pt, conference, english]{IEEEtran} 
\IEEEoverridecommandlockouts
\usepackage{graphicx}
\usepackage{epstopdf}
\usepackage{color}
\usepackage{epsfig}
\usepackage{subfigure}
\usepackage{amssymb,amsmath,amsfonts,layout,graphicx}
\usepackage{makeidx}
\usepackage{babel}
\usepackage{sublabel}
\usepackage{cases}
\usepackage{algorithm}
\usepackage{algorithmic}
\usepackage{multirow}
\usepackage{nomencl}

\usepackage
[
        letterpaper,
        left=1.91cm,
        right=1.91cm,
        top=1.91cm,
        bottom=2cm,
]
{geometry}

\newtheorem{problem}{Problem}

\newtheorem{proposition}{Proposition}

\newtheorem{definition}{Definition}

\newtheorem{remark}{Remark}
\newtheorem{example}{Example}

\DeclareMathOperator*{\argmax}{arg\,max}

\usepackage{cite}

 \title{A Mechanism Design Approach for Coordination of Thermostatically Controlled Loads}
  \author{Sen Li, \emph{Student Member}, \emph{IEEE}, Wei Zhang, \emph{Member}, \emph{IEEE}, Jianming Lian, \emph{Member}, \emph{IEEE}, and \\ Karanjit Kalsi, \emph{Member}, \emph{IEEE}
 \thanks{S. Li, and W. Zhang are with the Department of Electrical and Computer Engineering, Ohio State University, Columbus, OH 43210. Email: \{li.2886, zhang.491\}@osu.edu}
  \thanks{J. Lian and K. Kalsi are with the Electricity Infrastructure Group, Pacific Northwest National Laboratory, Richland, WA 99354. Email: \{jianming.lian, karanjit.Kalsi\}@pnnl.gov}
 } 
 
\begin{document} 
\maketitle
\begin{abstract}
This paper focuses on the coordination of a population of thermostatically controlled loads (TCLs) with unknown parameters to achieve group objectives. The problem involves designing the device bidding and market clearing strategies to motivate self-interested users to realize efficient energy allocation subject to a peak energy constraint. 
 This coordination problem is formulated as a mechanism design problem, and we propose a mechanism to implement the social choice function in dominant strategy equilibrium. The proposed mechanism consists of a novel bidding and clearing strategy that incorporates the internal dynamics of TCLs in the market mechanism design, and we show  it can realize the team optimal solution.
A learning scheme is proposed to address the unknown load model parameters. Numerical simulations are performed to validate the effectiveness of the proposed coordination framework.

\end{abstract}
\vspace*{0.1cm}

\section{Introduction}
Demand response has attracted considerable research attention in recent years, and is regarded as one of the most important means to improve the efficiency and reliability of the future smart grid. A natural way to achieve demand response is through various pricing schemes, such as Real Time Pricing (RTP), Time of Use (TOU) and Critical Peak Pricing (CPP) \cite{chao2010price}, \cite{allcott2009real}. Many validation projects \cite{faruqui2010impact} have been carried out to demonstrate the performance of these pricing schemes in terms of payment reduction, load shifting, and peak shaving.  These price-based methods either directly pass the wholesale energy price to end-users \cite{allcott2009real} or design pricing strategies in heuristic ways \cite{wolak2007residential}. It is thus hard to achieve predictable and reliable aggregated response, which is essential in various demand response applications, such as energy capping, load following, frequency regulation, among others.

To achieve accurate and reliable load response, aggregated load control has been extensively studied in the literature. A simple form of aggregated load control is the direct load control (DLC), where the aggregator can remotely control the operations of residential appliances based on the agreement between customers and the utility company. While traditional DLC is mainly concerned with peak load management \cite{hsu1991dispatch}, \cite{salehfar1991production}, recent research effort focuses more on the modeling and control of different kinds of aggregated loads, such as data center servers \cite{chen2013real,sen2014integrated},  hybrid electrical vehicles \cite{liu2013planning},\cite{han2010development} and thermostatically controlled loads  \cite{hao2013aggregate,TPS_13,lu2012evaluation,burke2008robust}, to participate in various demand response programs.  Some of these DLC methods require fast communications between the aggregator and individual loads. The communication overhead can be reduced using advanced state estimation algorithms \cite{mathieu2013state,vrettos2014demand} that can accurately estimate load state information without frequently collecting measurements from the loads.

Another important paradigm of aggregated load control is the market-based coordination. It borrows ideas from economics \cite{mas1995microeconomic} to coordinate a group of self-interested users to achieve desired aggregated load response \cite{fahrioglu2000designing}, \cite{samadi2012advanced}. Different from DLC, the market-based coordination affects the load response indirectly via an internal price signal. The internal price can be dramatically different from the wholesale price due to specific group objectives. For instance, in \cite{ygge1999decentralized} and \cite{ygge1997making}, a market-based approach is proposed  to efficiently allocate thermal resources among offices only based on local information. In \cite{kok2012dynamic} and \cite{kok2010intelligence}, a multi-agent based control framework is proposed to integrate distributed energy resources for various coordination objectives. A distributed algorithm is developed in \cite{chen2010two} and \cite{li2011optimal}  for the utility company and users to jointly determine optimal prices and demand schedules via an iterative bidding and clearing process. In \cite{paschalidis2012demand},  a group of smart buildings are coordinated through an internal price signal to provide frequency regulation services to the ancillary market. 
In addition, the Pacific Northwest National Laboratory  launched the GridWise$\textsuperscript{\textregistered}$ demonstration project to validate the market-based coordination strategies for residential loads\cite{fuller2011analysis}. The demonstration project involved 112 residential houses in Washington and Oregon, and showed that the market-based coordination strategies could reduce the utility demand and congestion at key times.

 Although the aggregated dynamics of TCLs may significantly affect the performance of the control strategies, many existing market-based coordination strategies either neglect this internal dynamics or use a simplified model to characterize it. In this paper, we consider the coordination of a group of TCLs to maximize the social welfare subject to a peak energy constraint, where the internal dynamics of TCLs are taking into account.
This coordination problem poses several challenges. First, the  user utilities are private information, making it rather challenging for the coordinator to achieve group objectives with incomplete information. Second, many existing works \cite{li2011optimal}, \cite{chen2012optimal} require multiple iterations between the agents and the coordinator to achieve the optimal social outcome. The real time implementation of such coordination algorithms requires considerable communication resources. Third, a lot of existing literature assumes accurate  load models with known parameters. However, the Gridwise$\textsuperscript{\textregistered}$ demonstration project \cite{fuller2011analysis} suggests this is not always the case. In practice, the information each user sends to the coordinator can only depend on local measurements, such as room temperature and ``on/off" state. Therefore, an estimation scheme is needed for the users to compute their bids only based on online measurements.

The key contribution of this paper lies in the development of a market-based coordination framework for residential air conditioning loads with a systematic consideration of all the aforementioned challenges. In this paper, we formulate the coordination problem as a mechanism design problem \cite{mas1995microeconomic,maskin2008mechanism}. The price-responsive loads are modeled as individual utility maximizers, while the group objective is encoded in the social choice function, which is to maximize the social welfare subject to a peak energy constraint. We propose a mechanism and show it can implement the social choice function in dominant strategy equilibrium. Such solution concept does not require iterative information exchanges between the coordinator and the individual loads, and can be implemented with limited communication resources. The proposed mechanism contains a novel bidding and clearing strategy that incorporates the internal dynamics of the TCLs into the market mechanism design, and we show that it can realize the team optimal solution. 

Different from many existing works, the problem is addressed with a systematic consideration of various practical factors, such as heterogeneous load dynamics, private information of individual users, unknown parameters of the load model, communication resources for the information exchange, etc. All these factors are brought up based on the observations in the GridWise$\textsuperscript{\textregistered}$ demonstration project \cite{fuller2011analysis}. They are important not only for customer privacy protection and the end user engagement, but also for the cost-effective implementation of the real-time control strategies.
Once our framework is properly implemented, it can accurately achieve the desired load responses, and improve the operational efficiency of the  distribution system in an economically feasible way.

\begin{figure*}[!ht]%
\begin{minipage}[b]{0.32\linewidth}
\centering
\includegraphics[width = 1\linewidth]{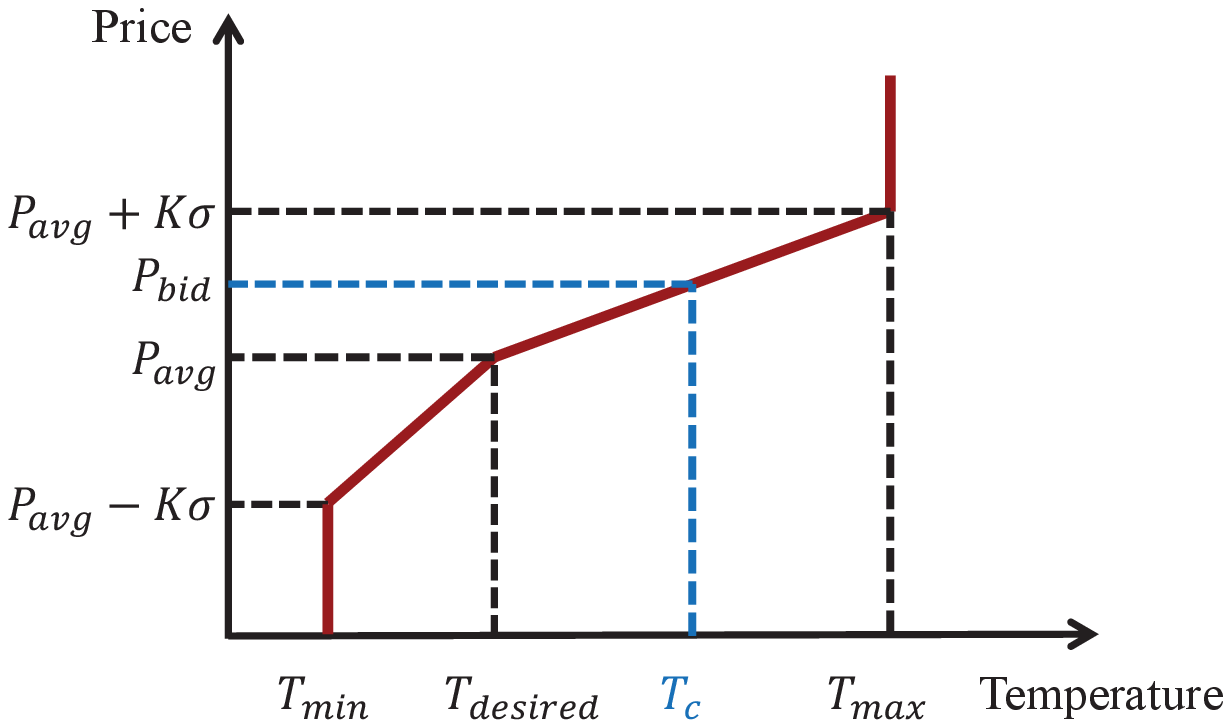}
\caption{The controller measures its current temperature $T_c$ and submits a bid $P_{bid}$ to the coordinator using this curve.}
\label{fig:bidding}
\end{minipage}
\begin{minipage}[b]{0.01\linewidth}
\hfill
\end{minipage}
\begin{minipage}[b]{0.32\linewidth}
\centering
\includegraphics[width = 1\linewidth]{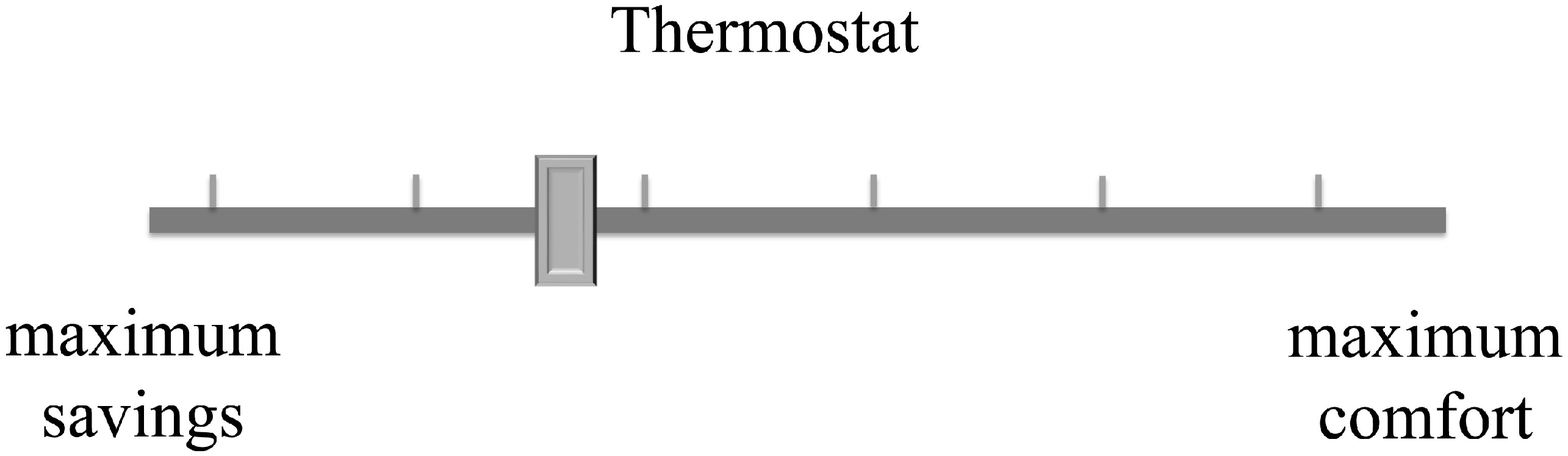}
\caption{User interface used in the GridWise$\textsuperscript{\textregistered}$ demonstration project \cite{fuller2011analysis}.}
\label{fig:slidingbar}
\end{minipage}
\begin{minipage}[b]{0.01\linewidth}
\hfill
\end{minipage}
\begin{minipage}[b]{0.32\linewidth}
\centering
\includegraphics[width = 0.8\linewidth]{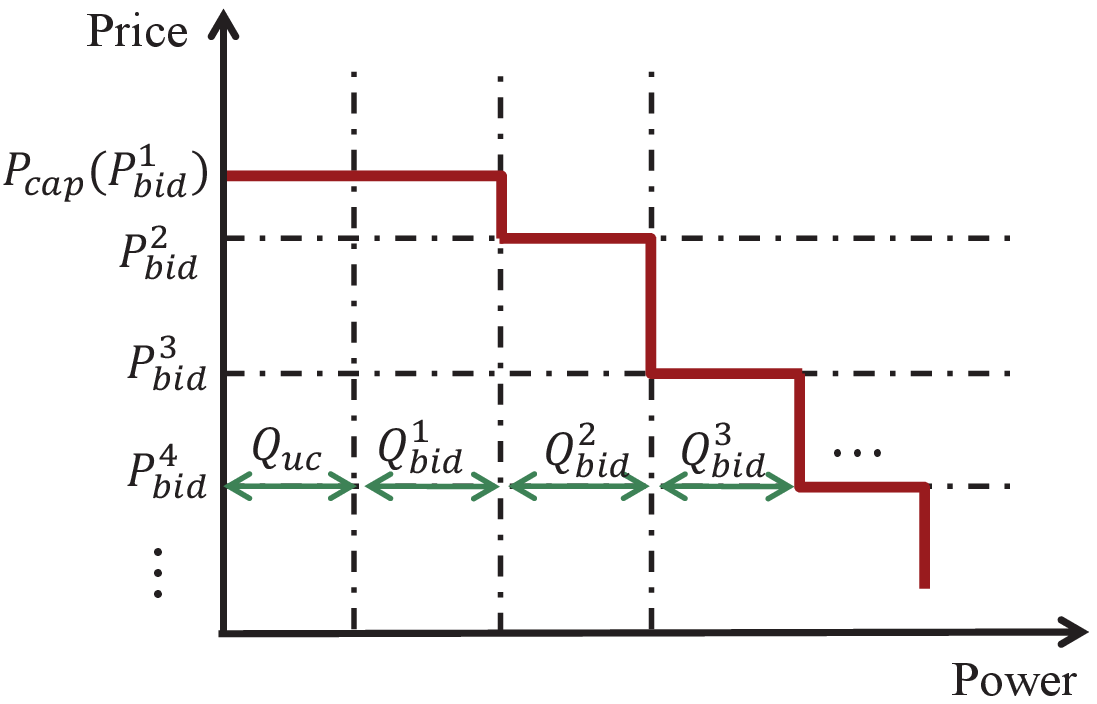}
\caption{The demand curve based on the user bids, where $P_{bid}^i$ is the bidding price sequence in decreasing order, and $Q_{bid}^i$ is the power of the most recent on cycle.}
\label{fig:demandcurve}
\end{minipage}
\end{figure*}

\begin{figure*}[!ht]%
\begin{minipage}[b]{0.32\linewidth}
\centering
\includegraphics[width = 0.8\linewidth]{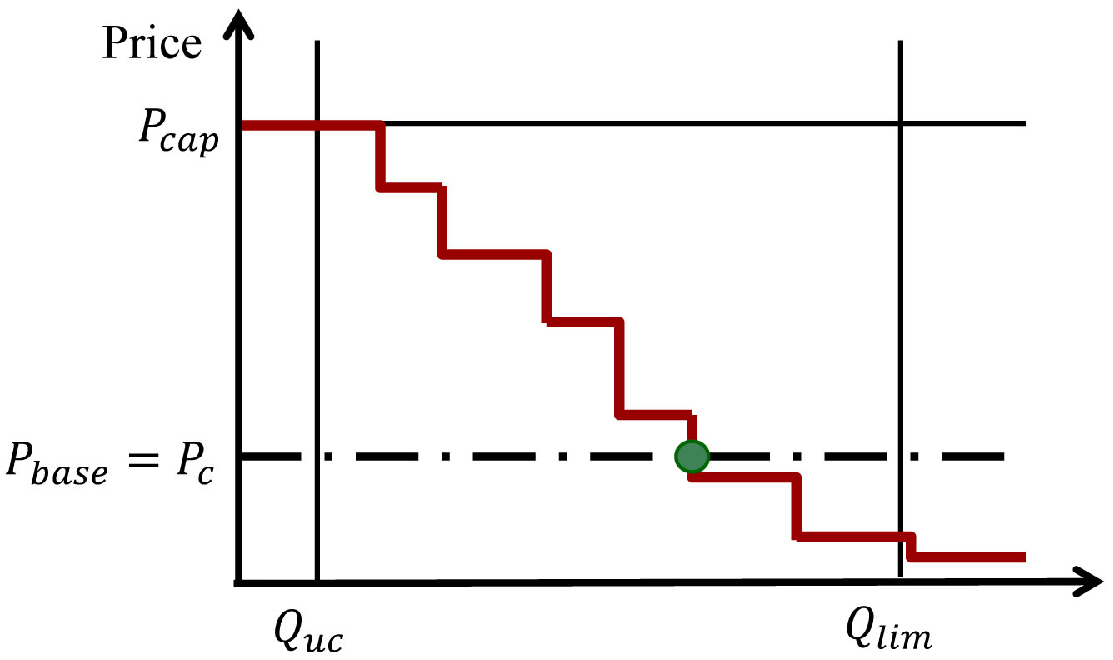}
\caption{The demand curve constructed based on all the bids. If the total demand is less than the feeder capacity constraint, then the clearing price is equal to the base price.}
\label{fig:demandcurve1}
\end{minipage}
\begin{minipage}[b]{0.01\linewidth}
\hfill
\end{minipage}
\begin{minipage}[b]{0.32\linewidth}
\centering
\includegraphics[width = 0.75\linewidth]{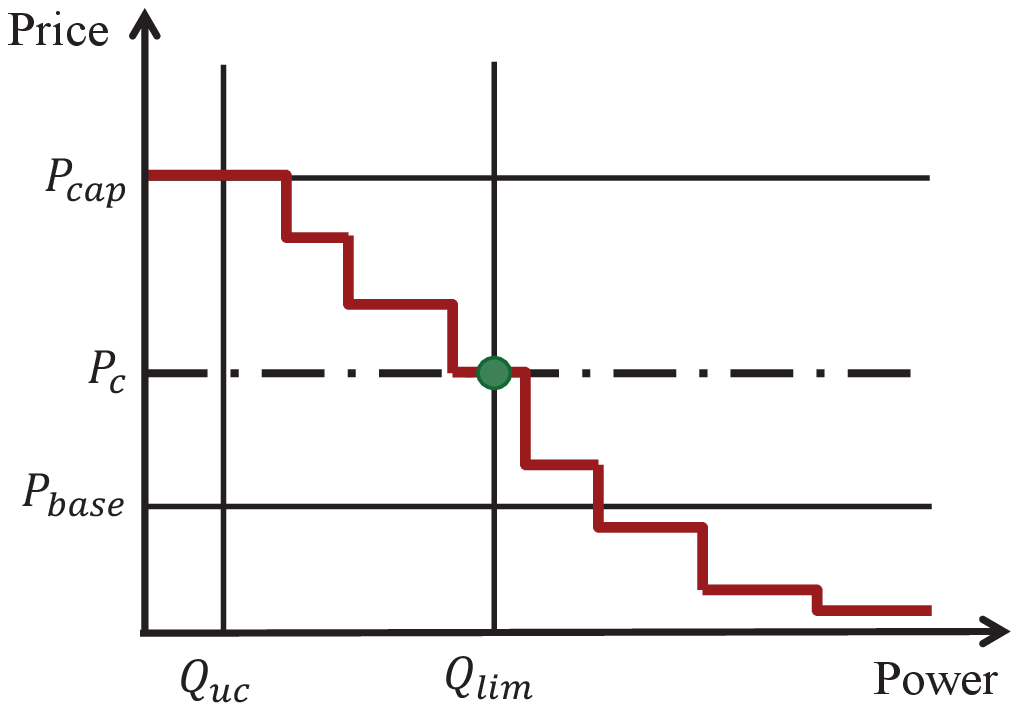}
\caption{When the total demand is greater than the feeder power constraint, then the clearing price is determined by the intersection of demand curve and feeder capacity constraint.}
\label{fig:demandcurve2}
\end{minipage}
\begin{minipage}[b]{0.01\linewidth}
\hfill
\end{minipage}
\begin{minipage}[b]{0.32\linewidth}
\centering
\includegraphics[width = 0.9\linewidth]{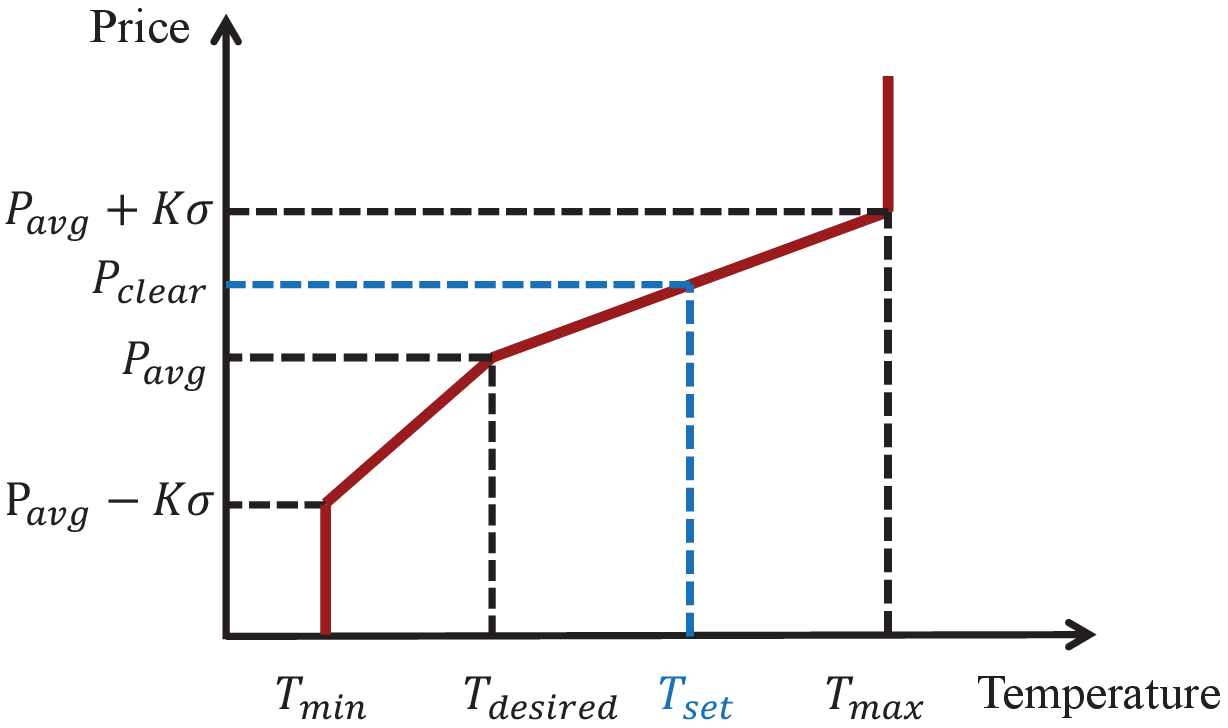}
\caption{The user response to the price. For any given price, the devices determine the temperature setpoint according to this curve. }
\label{fig:response}
\end{minipage}
\end{figure*}

The rest of the paper proceeds as follows. A motivating example based on a real-world demonstration project is presented in Section II, followed by a problem formulation in Section III. A mechanism is constructed in Section IV to implement the optimal energy allocation. A joint state-parameter estimation framework is presented in Section V, followed by simulations results in Section VI and some concluding remarks in Section VII.

\section{Motivating Example}

The framework proposed in this paper is largely motivated by the Pacific Northwest GridWise$\textsuperscript{\textregistered}$ demonstration project \cite{fuller2011analysis}, where  a 5-minute double-auction market is created to coordinate  a group of TCLs to cap the aggregated peak energy.  
Each device is equipped with a smart thermostat that can measure the room temperature and communicate with the coordinator. Before each market period, the device measures its room temperature, $T_c$, and submits a bid to the coordinator. The bid should consist of the load power and the bidding price. Since the rated power of the load is different from its actual power due to environmental disturbances, in practice each device is required to bid the measured average power of the most recent market period during which the load is on. The bidding price is determined by a bidding curve shown in Fig. \ref{fig:bidding}, where $P_{avg}$ is the average clearing price of certain price history (e.g., 24 hours), $\sigma$ is the standard variation of the clearing prices during the given history, and $T_{min}$, $T_{desired}$ and $T_{max}$ are user-specified minimum, desired, and maximum temperature, respectively. We denote the bidding power and price as $Q_{bid}$ and $P_{bid}$, respectively.  
In addition, each user can specify energy use preferences through a smart thermostat interface (see Fig. \ref{fig:slidingbar}). This user preference will affect the slope of the bidding curve.

The coordinator collects all the bids and orders the bids in a decreasing sequence, $P_{bid}^1,\ldots, P_{bid}^N$. With the associated power sequence, $Q_{bid}^1,\ldots,Q_{bid}^N$, a demand curve can be constructed to map the clearing price to aggregated power. Fig. \ref{fig:demandcurve} illustrates how the demand curve is constructed. This curve is then used to determine the market clearing price that respects the feeder capacity constraint: when the total demand is less than the feeder capacity, the market clearing price is equal to the base price, $P_{base}$ (Fig. \ref{fig:demandcurve1}), which is the wholesale energy price plus a retail modifier as defined by the tariff of American Electric Power (AEP) \cite{aeptariff}; otherwise  the market price, $P_c$, is determined by the intersection of the demand curve and the feeder capacity constraint (Fig. \ref{fig:demandcurve2}). 

After the market is cleared, each device receives the energy price and adjusts its setpoint, $T_{set}$, according to a response curve as shown in Fig. \ref{fig:response}. This setpoint modifies the system dynamics and affects the temperature trace of the TCL, and therefore affects the bid of each user for the next market period. Notice that all the bidding and user response processes are executed by a programmable controller, and the user only needs to specify his/her preferences via the thermostat interface. To initialize the market process, the user needs to specify $T_{min}$, $T_{max}$, $T_{desired}$ and $K$, the device needs to measure the temperature and the power of the last ``on" cycle, and the coordinator needs to collect all the bids, estimate the power of the unresponsive loads, $Q_{uc}$, and the feeder capacity constraint, $D$.  

Apart from the GridWise$\textsuperscript{\textregistered}$ project, a similar demonstration project is also implemented in AEP, Ohio \cite{AEPdemo}, which involves more households and more sophisticated market bidding design. These projects provide insights for the coordination of residential loads from the practical point of view. However, the bidding and pricing strategies are designed in a heuristic way, which may result in constraint violations and market inefficiencies. To address these challenges, there is a strong need to develop a general coordination framework that can serve as a theoretical foundation to improve the performance of the control scheme and help to design other similar market-based coordination strategies.






\section{Problem Formulation}
Consider a coordination problem for a group of TCLs, where the coordinator allocates energy to users to maximize the social welfare subject to a total energy constraint. Each device is assumed to be equipped with a smart thermostat that has two main functions. First, it allows the user to specify energy use preferences via an interface such as the sliding bar shown in Fig. \ref{fig:slidingbar} to indicate one's trade-off between comfort and cost. Second, before each market period it submits a bid  to the coordinator based on user's preference and local device measurement, such as power consumption, ``on/off" states, and local temperature. 
The coordinator collects the user bids, determines the energy price, and broadcasts the price to all the devices. Each device will then adjust the temperature setpoint in response to the energy price to maximize the individual utility. This will modify the system dynamics and therefore affect the user bids for the next period. In the considered scenario, we assume that each user is a price taker, namely, an individual user's decision will not significantly affect the market price. This is a standard assumption when the market involves a large number of players \cite[chap. 12.F]{mas1995microeconomic}, \cite{conejo2002price}, \cite{de2007strategic}.

The rest of this section provides formal mathematical descriptions of the main components of the proposed framework.

\subsection{User Preferences and Utility}
Assume that there are $N$ self-interested users. Each user needs to determine the temperature setpoint to obtain an energy allocation that maximizes his individual utility (the user's comfort minus the electricity cost). In other words, each user is confronted with the trade-off between comfort and electricity cost: when the electricity price is high, the device will adjust the temperature setpoint to save electricity cost at the sacrifice of user comfort. Formally,  a function $V_i:\mathit{\mathbb{R}\rightarrow\mathbb{R}}$ can be used to represent the comfort level for each user with energy allocation $a_i$. Assume that $V_i(a_i)$ is concave, continuously differentiable, $V_i(0)=0$ and $V'_i(0)>0$. Let $\theta_i(t_k)$ represent the private information of user $i$. Denote $E_i^m$ as the energy consumption for the $i$th load if it is ``on" during the entire period, which gives $a_i\leq E_i^m$. The individual utility maximization problem can be formulated as follows:
\begin{align}
\label{eq:individualopt}
&\max_{a_{i}} \quad V_i(a_{i};\theta_i(t_k))-P_ca_{i}\\
&\ \text{subject to: }  0\leq a_{i}\leq E_{i}^{m} ,  \nonumber
\end{align}
where $P_c$ is the energy price. Let $h_i:\mathit{\mathbb{R}\rightarrow\mathbb{R}}$ be the optimal solution to the optimization problem (\ref{eq:individualopt}), we have:
\begin{equation}
\label{eq:userresponse}
h_{i}(P_c;\theta_i(t_k))=\argmax_{0\leq a_{i} \leq E_i^{m}} V_i(a_{i};\theta_i(t_k))-P_ca_{i} .
\end{equation}
We assume that $h_i$ is continuous and non-increasing with respect to $P_c$ for each $i=1,\ldots,N$.   
Notice that the user can not directly choose his optimal energy allocation. Instead, he can only determine the temperature setpoint, which affects the energy consumption through the load dynamics.

\subsection{Individual Load Dynamics}
Let $\eta_i(t)\in \mathbb{R}^{n}$ be the continuous  state of the $i$th load. Denote $q_i(t)$ as the ``on/off" state: $q_i(t)=0$ when the TCL is off, and $q_i(t)=1$ when it is on.
For both ``on'' and ``off'' states, the thermal dynamics of a TCL system can be typically modeled as a linear system:
\begin{align}
\dot{\eta}_i(t)=\begin{cases}
A_{i}\eta_{i}(t)+B_{on}^{i}  & \text{ if } q_i(t)=1\\
A_{i}\eta_{i}(t)+B_{off}^{i}  & \text{ if } q_i(t)=0 .
\end{cases}
\label{eq:inddynamics}
\end{align}
Many existing works use a first-order linear system to capture the TCL dynamics  \cite{hao2013aggregate}, \cite{mathieu2013state}, \cite{vrettos2014demand}, where $\eta_i(t)$ only consists of the room temperature.
Although the first-order model is adequate for small TCLs such as refrigerators, it is not appropriate for residential air conditioning systems, which require a 2-dimensional linear system model incorporating both air and mass temperature dynamics \cite{TPS_13}.
Such a second-order model is typically referred to as the Equivalent Thermal Parameter (ETP) model \cite{sonderegger1978dynamic}.
In this paper we focus on the second-order ETP model, which includes the first-order model as a special case. Let $\varphi_i=[A_i,B_{on}^i,B_{off}^i]^T$ be the model parameters. Typical values of these parameters and the factors that affect these parameters can be found in \cite{TPS_13}.

\begin{figure}[t]
\centering
\includegraphics[width =0.85\linewidth]{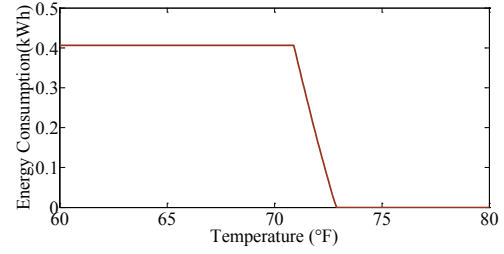}
\caption{Energy consumption of the TCL during a market clearing cycle as a function of the temperature setpoint.}
\label{energyfunction}
\end{figure}

The power state of the TCL is typically regulated by a hysteretic controller based on the control deadband $[u_i(t)-\delta/2,u_i(t)+\delta/2]$, where $u_i(t)$ is the temperature setpoint of the $i$th TCL and $\delta$ is the deadband. Let $T_c^i(t)$ denote the room temperature of the $i$th load.
In the cooling mode, the load is turned off  when $T_c^i(t)\leq u_i(t)-\delta/2$, and it is turned on when $T_c^i(t)\geq u_i(t)+\delta/2$, and remains the same power state otherwise. This hysteretic control policy can be described as:
\begin{align}
q_i(t^+) = \begin{cases}
1 & \text{ if } T_c^i(t)\geq u_i(t)+\delta/2\\
0 & \text{ if } T_c^i(t)\leq u_i(t)-\delta/2\\
q_i(t)  & \text{ otherwise } .
\end{cases}
\label{eq:hvaccontrol}
\end{align}
For notation convenience, we define a hybrid state $z_i(t)=[\eta_i(t),q_i(t)]^T$, which consists of both the temperature and the ``on/off" state of the load. Let $[t_k,t_k+T]$ be the $k$th market period, then the energy consumption of each load during the $k$th period depends on the system state and setpoint control $u_i(t)$. In this case, the private information consists of system state and model parameters. Therefore, the energy consumption of each load can be represented as $e_i(u_i(t_k),z_i(t_k),\varphi_i)$. This energy consumption function can be derived by calculating the portion of time that the system is on over the entire market period (details of this calculation  are presented in Section IV). An example is shown in Fig \ref{energyfunction}, where a second-order ETP model is used and the initial room temperature is $72.8^\circ$F.  Let $\theta_i(t_k)=(z_i(t_k),\varphi_i)$ be the overall private information of load $i$, then the energy function can be written as $e_i(u_i(t_k),\theta_i(t_k))$. Notice that the private information for users is time varying, as it contains the system state. 

\begin{figure}[t]
\centering
\includegraphics[width = 0.85\linewidth]{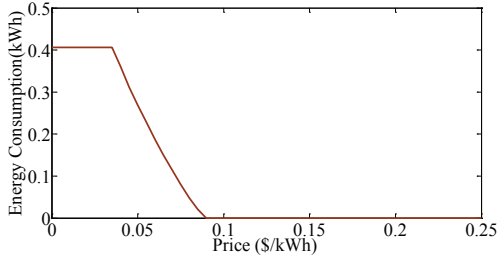}
\caption{Energy consumption of the TCL during a market clearing cycle as a function of the energy price.}
\label{energytprice}
\end{figure}

After the market is cleared, each user wants to determine the control action $u_i(t_k)$ such that the resulting energy consumption equals the optimal solution to (\ref{eq:individualopt}). Since the optimal control depends on the energy price, we can define a user response function, $\Lambda_i:\mathit{\mathbb{R}\rightarrow\mathbb{R}}$ with $u_i(t_k)=\Lambda_i(P_c)$. Therefore, the optimal energy allocation function $h_i$ as defined in (\ref{eq:userresponse}) should satisfy the following:
\begin{equation}
\label{eq:interp}
 h_i(\cdot;\theta_i(t_k))=e_i(\Lambda_i(\cdot),\theta_i(t_k)) .
\end{equation} 
The left-hand side of equation (\ref{eq:interp}) represents the optimal energy allocation for a given price, while the right-hand  side arises from the physical property of the individual loads, and indicates that the user can specify the control action $u_i$ to match  the actual energy consumption to the optimal allocation. An example of function $h_i$ is shown in Fig. \ref{energytprice}, where the response curve is piecewise linear (as shown in Fig. \ref{fig:bidding}) and the initial room temperature is $72.8^\circ$F. To derive the function $h_i(\cdot;\theta_i(t_k))$, we first determine the control setpoint based on the market price using the response curve (Fig. \ref{fig:bidding}), then calculate the corresponding energy consumption based on the energy function $e_i(\cdot,\theta_i(t_k))$. Since the energy function $e_i(\cdot,\theta_i(t_k))$ depends on the system dynamics (\ref{eq:inddynamics}) and the control policy (\ref{eq:hvaccontrol}), the load dynamics are incorporated in function $h_i$ through this process.

\subsection{Problem Statement}
The coordinator obtains energy from the wholesale market at a cost denoted as $C\left(\sum_{i=1}^N a_i\right)$. We assume that $C(\cdot)$ is differentiable and convex. The energy is then allocated to users via a price signal to maximize the social welfare, which can be defined as $\sum_{i=1}^{N}V_{i}(a_{i};\theta_i(t_k))-C(\sum_{i=1}^{N}a_{i})$. Therefore, the coordinator's optimization problem can be formulated as follows:
\begin{align}
\label{eq:coordinatoropt}
&\max_{a_1,\ldots,a_N} \sum_{i=1}^{N}V_{i}(a_{i};\theta_i(t_k))-C\left(\sum_{i=1}^{N}a_{i}\right) \\
&\text{subject to: } 
\begin{cases}
\sum_{i=1}^{N}a_{i}\leq D \\ \nonumber
0\leq a_{i}\leq E_{i}^{m}, \forall i=1,\ldots,N \\ \nonumber
a_i=h_i(P_c;\theta_i(t_k)), \forall i=1,\ldots,N,
\end{cases}
\end{align}
where $D$ is the maximum energy for the aggregated loads. Without loss of generality, we assume that $D\leq N E_i^m$. Note that the feeder capacity constraints considered in the GridWise$\textsuperscript{\textregistered}$ demonstration project can be represented by the total energy constraint. This is because the feeder capacity constraint is mainly due to the consideration of the thermal characteristics of the feeder. The instantaneous power can exceed the feeder power limit without causing damages to the grid, as long as the energy over a certain period is effectively capped to protect the feeder from overheating. 

The optimization problem (\ref{eq:coordinatoropt}) defines a Stackelberg game \cite{basar1995dynamic}, where the coordinator first makes control decisions to maximize the social welfare, then the individual users choose energy consumption to maximize individual utility based on the coordinator's control decisions. In such Stackelberg games, the upper bound on the social welfare can be typically characterized by the team optimal solution \cite{basar1995dynamic}, which is the optimal solution to the following team problem:
\begin{align}
\label{eq:teamopt}
&\max_{a_1,\ldots,a_N} \sum_{i=1}^{N}V_{i}(a_{i};\theta_i(t_k))-C\left(\sum_{i=1}^{N}a_{i}\right) \\
&\text{subject to: } 
\begin{cases}
\sum_{i=1}^{N}a_{i}\leq D \\ \nonumber
0\leq a_{i}\leq E_{i}^{m}, \forall i=1,\ldots,N ,
\end{cases}
\end{align}
In the above team problem, the coordinator and the users cooperatively maximize the social welfare subject to the peak energy constraint. In general, the team solution results in a higher social welfare than the solution to (\ref{eq:coordinatoropt}), since the coordinator's optimization problem (\ref{eq:coordinatoropt}) is more restrictive: one only needs to find an energy allocation to maximize the social welfare to solve the team problem, while in the coordinator's optimization problem, we also need to find a price to satisfy the additional constraint in (\ref{eq:coordinatoropt}). However,  such a clearing price may not always exist for an arbitrarily given team optimal solution.

\begin{example}
As an example, consider two users with $V_1(a_1;\theta_1(t_k))=a_1$, $V_2(a_2;\theta_2(t_k))=3a_2$.  The energy cost for the coordinator is $C(a_1+a_2)=2a_1+2a_2$. The team problem is to maximize the social welfare subject to an energy constraint, i.e.:
\begin{align}
\label{eq:globalopt1}
&\max_{a_1,a_2} \sum_{i=1}^{2}V_{i}(a_{i};\theta_i(t_k))-C(a_1+a_2) \\
&\text{subject to: } 
\begin{cases}
a_1+a_2 \leq 1 \\ \nonumber
0\leq a_{i}\leq 2, \text{ for }  i=1,2
\end{cases}
\end{align}
The team optimal solution is $a_1=0$, $a_2=1$. However, according to (\ref{eq:individualopt}), given any energy price, $a_i$ is either 0 or 2. Therefore, the coordinator can not find a price to realize the team optimal solution.
\end{example}
\vspace{0.05in}

To address this concern, we introduce the concept of realizable energy allocation:
\begin{definition}
The energy allocation vector, $a=(a_1,\ldots,a_N)$, can be realized by $P_c$, if $a_i=h_i(P_c;\theta_i(t_k))$ for all $i=1,\ldots,N$.
\end{definition}

It is clear that not all the energy allocations can be realized. In this paper, 
we have assumed that $V_i$ is concave and continuously differentiable, and $h_i$ is continuous and non-increasing. We will show in Section V that under these conditions, there is always a price to realize the team optimal solution. In other words, the upper bound given by the team optimal solution is tight. Therefore, the problem of the paper can be formulated as follows:
\begin{problem}
Design the bidding and clearing strategy, such that the cleared price realizes the team optimal solution $a^*$.
\end{problem}

The coordinator's optimization problem (\ref{eq:coordinatoropt}) can not be directly addressed using standard optimization techniques, since the individual valuations are unknown to the coordinator. For this reason, to achieve the group objectives, the coordinator needs to design a bidding strategy to collect information from the individual users, and then determine the price based on the user bids. 

\begin{remark}
The market design for many traditional assets is well-understood. For instance, in energy market,  generators can be simply characterized by an output range depending on its ramp rate during each market period. However, the internal dynamics of TCLs are more complex and depend more on the environment, and thus cannot be handled in the same way. Therefore, an important contribution of this paper is to incorporate the dynamics of TCLs in the energy market design.
In addition, although this paper only considers the load dynamics within one market period, it is an important step towards establishing a fully dynamic version of the problem where multiple market periods are taken into account. 
\end{remark}

\section{A Mechanism Design Framework}
In this section, we adopt the mechanism design approach to solve Problem 1. First the problem is formulated as a mechanism design problem, then a mechanism is constructed to implement the desired social outcome. In addition, a realistic bidding strategy with a simplified message space is proposed to reduce the communication overhead.

\subsection{The Mechanism Design Problem}
Mechanism design studies how to aggregate the individual preferences into  a social choice while the individual's actual preferences are not publicly observable. In a mechanism design problem, each user is assumed to selfishly take actions to maximize the individual utility, while the coordinator makes the collective choice that achieves various group objectives. Since the individual utility is unknown to the coordinator, he can require each user to submit a bid to collect information. In this case, the key problem for the coordinator is to align individual objectives with system-level objectives. In other words, a proper bidding and pricing strategy needs to be designed, such that when each user selfishly maximizes the individual utility, the resulting outcome also achieves the desired group objectives (for example, maximizes the social welfare). The rest of this subsection introduces basic concepts in mechanism design. 

Let $x\in X$ be the outcome of the mechanism that consists of the energy allocation and the energy price, i.e., $x=(a_1,\ldots,a_N,P_c)$. The utility of each  user (comfort minus electricity cost) depends on the outcome. Moreover, we assume that at time $t_k$, each user can privately observe his utility, $U_i$, over different outcomes. In other words, we can model this by supposing that user $i$ privately observes a parameter $\theta_i$ that determines his utility. Notice that we drop the dependence of $\theta_i$ on $t_k$ throughout the rest of the paper for notation convenience. In mechanism design, $\theta_i\in \Theta_i$ is usually referred to as the user $i$'s {\em type} \cite[p. 858]{mas1995microeconomic}, where $\Theta_i$ denotes the set of all the possible {\em types}. In our problem, the user {\em type} contains the system state, $z_i(t_k)$, and the model parameter, $\varphi_i$, in particular:
\begin{equation}
\label{eq:utility}
 U_i(x;\theta_i)=V_i(a_i;\theta_i)-P_ca_i ,
\end{equation} 
where $\theta_i=[z_i(t_k),\varphi_i]$.

As the user preferences are private, to determine the optimal energy price, the coordinator also needs to require each user to submit a bid to reveal some information. Formally, this can be formulated as a message space  $M=M_1\times \cdots \times M_N$, where $M_i$ denotes the space of messages (bids) the $i$th user can communicate to the coordinator.  The structure of $M_i$ depends on particular applications. For example, in the demonstration project, each device submits a price and a quantity, then we have $(P_{bid}^i,Q_{bid}^i)\in M_i$. In \cite{chen2010two} each device submits the slope of the demand curve, $\beta_i$, in which case $\beta_i\in M_i$. After collecting all the user bids, the market is cleared with an energy price and a corresponding energy allocation. The clearing strategy can be represented by an outcome function, $g:\mathit{M\rightarrow X}$, that maps the user bids to an outcome, $x$. 
The message space and the outcome function together fully characterize the rules governing the procedure for making the collective choice. This is typically referred to as a {\em mechanism}\cite{mas1995microeconomic}, which can be denoted as $\Gamma=(M_1,\ldots,M_N,g(\cdot))$.

Each user observes $\theta_i$ privately and determines what to bid to maximize his utility. This process can be represented by a bidding strategy $m_i: \Theta_i\rightarrow M_i$ that maps the user {\em type} to a message. There are many solution concepts for a mechanism, such as Nash equilibrium, Bayesian Nash equilibrium, etc. Of particular interest to our framework in this paper is the dominant strategy equilibrium. Denote $m_{-i}$ as the collection of strategies of all the users other than $i$, then the dominant strategy equilibrium is defined as follows:
\begin{definition}[\text{Dominant Strategy Equilibrium} \cite{mas1995microeconomic}]
The strategy profile $(m_1^*(\cdot),\ldots,m_N^*(\cdot))$ is a dominant strategy equilibrium of mechanism $\Gamma=(M_1,\ldots,M_N,g(\cdot))$ if for all $i$ and all $\theta_i\in \Theta_i$, $U_i(g(m_i^*(\theta_i),m_{-i}),\theta_i)\geq U_i(g(m'_i(\theta_i),m_{-i}),\theta_i)$ for all $m'_i(\theta_i)\in M_i$ and all $m_{-i}\in M_{-i}$.
\label{def:dominantstrategy}
\end{definition}
\vspace{0.05in}

The equilibrium strategy characterizes the individual's self-interested behavior: each user is an individual welfare maximizer. However, in the coordinator's point of view, a more interesting question is to find the best choice for the overall social welfare. For this reason, a social choice function $f:\mathit{\Theta\rightarrow X}$ can be defined to represent the desired social outcome of the coordinator. More specifically, $f(\cdot)$ determines what outcome will be chosen by the coordinator when he knows all the private information. In our problem, $f$ consists of  the optimal price to the optimization problem (\ref{eq:coordinatoropt}) and the resulting energy allocation.  If we define $\theta=(\theta_1,\ldots,\theta_N)$, the conflict between the personal interest and social interest can be captured by the concept of {\em implementation}:
\begin{definition}[\text{Implementation} \cite{mas1995microeconomic}]
A mechanism $\Gamma=(M_1,\ldots,M_N,g(\cdot))$ implements the social choice function $f(\cdot)$ in dominant strategies if there exists a dominant strategy equilibrium $m^*(\cdot)$ of $\Gamma$, such that $g(m_1^*(\theta_1),\ldots,m_N^*(\theta_N))=f(\theta)$ for all $\theta\in \Theta$.
\end{definition}
\vspace{0.05in}

In the above definition, $g(m_1^*(\theta_1),\ldots,m_N^*(\theta_N))$ represents the resulting outcome of individual maximization, while $f(\theta)$ denotes the desired social outcome. The concept of implementation characterizes the social choice that can be realized when all the users take actions to selfishly maximize the individual utility. To this end, Problem 1 can be equivalently stated as follows:
\begin{problem}
Design a mechanism to implement the social choice function $f(\cdot)$ that maximizes the social welfare subject to a peak energy constraint, i.e., $f(\theta)=(h_1(P_c^*;\theta_i(t_k)),\ldots,h_N(P_c^*;\theta_i(t_k)),P_c^*)$ and $P_c^*$ is the solution to the optimization problem (\ref{eq:coordinatoropt}). Furthermore, $P_c^*$ realizes the team optimal solution.
\end{problem}
\vspace{0.05in}

In the above mechanism design problem, the coordinator needs to design the message space and the market clearing rule such that the optimal social welfare can be implemented when each user selfishly maximizes the individual utility. In the meanwhile, the peak energy constraint needs to be respected.

\subsection{Constructing the Mechanism}
Let $f(\theta)=(a_1^*,\ldots,a_N^*,P_c^*)$ be the social choice function that maximizes the social welfare subject to the peak energy constraint. Specifically, $P_c^*$ is the optimal solution to (\ref{eq:coordinatoropt}), and $f(\theta)$ satisfies the following condition:
\begin{equation}
a_i^*=h_i(P_c^*;\theta_i), \forall i=1,\ldots,N .
\label{socialchoice}
\end{equation}
This subsection constructs a mechanism to implement $f(\cdot)$.
Consider a mechanism $\Gamma^*$, where each device is asked to submit function $h_i(\cdot;\theta_i)$. Since we have assumed that $h_i(P_c;\theta_i)$ is continuous and non-increasing with respect to $P_c$, the message space is the function space of all non-increasing and continuous functions. Notice that the user's actual bids may deviate from function $h_i$, unless they are motivated to bid $h_i$. Let $b_i(\cdot;\theta_i)$ be a non-increasing and continuous function that represents the user's actual bid.  The aggregated demand curve $b(\cdot;\theta)$ can be obtained by adding individual bidding functions, i.e., $b(\cdot;\theta)=\sum_{i=1}^N b_i(\cdot;\theta_i)$.  In this mechanism, each user is required to submit a function, which requires considerable communication resources. This bidding strategy will be simplified in the next subsection to reduce the communication overhead. 

Here we propose the following outcome function $g(b_1,\ldots,b_N)=(a_1^*,\ldots,a_N^*,P_c^*)$ to clear the market:
\begin{numcases}{}
a_i^*=b_i(P_c^*;\theta_i) \text{ for all } i=1,\ldots,N  \label{EQ_1} \\
P_c^*=\text{max} \big \{\bar{P},P^* \big \} \label{EQ_2}\\
P^*=C'\left( \sum_{i=1}\nolimits ^N a_i^*\right) \label{EQ_4}\\
b(\bar{P}, \theta)=D \label{EQ_3} ,
\end{numcases}
where $C'$ represents the derivative of the cost function $C(\cdot)$. According to (\ref{EQ_4}) and (\ref{EQ_3}), $P^*$ is the marginal production cost of procuring $\sum_{i=1}^N a_i^*$ amount of energy, while $\bar{P}$ is the energy price at which the aggregated demand is equal to the maximum allowed amount. Since $b_i$ is continuous and non-increasing, and we have assumed that $D\leq N E_i^m$, $\bar{P}$ exists. Intuitively, the social welfare is maximized when the market price equals the marginal production cost, i,e, $P_c^*=P^*$. However, in equation (\ref{EQ_3}), the function $b$ is non-increasing with respect to price, indicating that any feasible price that respects the feeder capacity constraint should be greater than $\bar{P}$. Therefore, in the proposed outcome function, the clearing price equals to $P^*$ whenever $P^*>\bar{P}$, and equals to $\bar{P}$ otherwise. When the energy price is determined, the allocation exactly follows the user bids, i.e., $a_i^*=b_i(P_c;\theta_i)$. For illustrating purpose, we construct the following example to show how to derive the optimal solution from the proposed clearing strategy. 
\begin{example}
Consider 100 users with $V_i=-\frac{1}{2}a_i^2+(i-P_c)a_i$. Assume that after proper scaling, the maximum energy consumption for each user is 1. The individual utility maximization problem can be formulated as follows:
\begin{align}
\label{eq:indopt}
&\max_{a_i} -\frac{1}{2}a_i^2+(i-P_c)a_i \\
&\ \text{subject to: } 0\leq a_{i}\leq 1 \nonumber\
\end{align}
The optimal solution to this problem is:
\begin{align}
a_i^* = \begin{cases}
0 & \text{ if } P_c\geq i\\
1 & \text{ if } P_c\leq i-1\\
i-P_c  & \text{ otherwise } .
\end{cases}
\label{eq:optimalsolution}
\end{align}
In addition, let us assume that the real time price is 20, and the maximum  5-minute energy due to the feeder capacity constraint is 50, i.e., $P^*=20$ and $D=50$.
According to (\ref{eq:optimalsolution}), when $P_c=99$, only the $100$th user consumes 1 unit of energy, and the aggregated energy is 1. When $P_c=98$, the $99$th and the $100$th user consume 1 unit of energy, respectively, and the corresponding aggregated energy is 2, and so forth.  Therefore, the price that corresponds to the energy limit is 50, i.e., $\bar{P}=50$. Since $\bar{P}>P^*$, we conclude that $P_c^*=\bar{P}$. 
\end{example}
\vspace{0.05in}

The rest of this subsection discusses some properties of the proposed mechanism.

\begin{proposition}
When each user is a price taker, the strategy profile $(h_1(\cdot;\theta_1),\ldots,h_N(\cdot;\theta_N))$ is a dominant strategy equilibrium of the proposed mechanism $\Gamma^*$.
\label{prop1}
\end{proposition}
\vspace{0.05in}
This result follows easily from the price taker assumption. Its proof can be found in the Appendix section. In the proposed mechanism, the optimal bid of each user does not depend on the bidding decisions of others. This is a very important property, since in our particular problem, each user does not know other user's preferences or actions. Therefore, if the bidding decision of one user has to depend on the action of another, then the equilibrium strategy can not be achieved unless all the users have accurate predictions on other user's action, which may not be a reasonable assumption. In addition, we also want to comment that Proposition 1 only holds when there are many users such that the influence of an individual user on the market price is negligible. In other cases (such as the oligopolistic market), the mechanism needs to be designed differently.

Now we can establish the following key property of the proposed mechanism:
\begin{proposition}
The proposed mechanism $\Gamma^*$ implements the social choice function $f(\cdot)$. Furthermore, the resulting market clearing price realizes the team optimal solution.
\label{prop2}
\end{proposition}
\vspace{0.05in}
The proof of this proposition can be found in Appendix.

\subsection{Realistic Bidding Strategy}
The proposed mechanism provides a general solution to the coordination problem formulated in this paper. In real-world applications, directly submitting function $h_i$  requires considerable communication resources, and might impinge on the customer privacy. Therefore, in this subsection we explore the structure of function $e_i(\cdot;\theta_i)$ and $h_i(\cdot;\theta_i)$ to simplify the message space and reduce the communication overhead.

In this paper we assume that the TCL consumes a constant power when it is ``on'', and consumes no energy when it is ``off''. For this reason, the energy consumption function $e_i(\cdot,\theta_i(t_k))$ can be derived by calculating the portion of time that the system is on during the entire market period. For example, assume that the system is ``on'' at the end of the $(k-1)$th period. When the initial temperature $\eta_i(t_k)$ is given, the state trajectory of the linear dynamic model (\ref{eq:inddynamics}) can be derived as $\eta_i(t)=e^{A_it}\eta_i(t_k)+A_i^{-1}(e^{A_it}-I)B_{on}$, where $\eta_i(t_k)=[\eta^{(1)}_i(t_k),\eta^{(2)}_i(t_k)]^T$, $\eta^{(1)}_i(t_k)=T_c^i(t_k)$ and $I$ is the identity matrix. When the trajectory hits the boundary of the control deadband defined in (\ref{eq:hvaccontrol}), the power state will switch and the system is off. Therefore, the trajectory of the system state $\eta_i(t)$ and the power state $q_i(t)$ for the entire period  can be derived, and the portion of time that the system is ``on'' can be  calculated based on $q_i(t)$.  In particular, consider a system in cooling mode. If the load is ``on'' at the end of the $(k-1)$th period, i.e., $q_i(t_k^-)=1$, we have the following (the case for $q_i(t_k^-)=0$ can be derived similarly):
\begin{align}
e_i(u_i(t_k),\theta_i(t_k)) = \begin{cases}
E_i^m & \text{ if } u_i(t_k)\leq T_f^i(t_k)+\delta/2\\
0 & \text{ if } u_i(t_k)\geq T_c^i(t_k)+\delta/2 \nonumber \\
E_i^m\times\alpha  & \text{ otherwise }, 
\end{cases}
\end{align}
where $\alpha=\int_0^T q_i(t)dt=\frac{T'}{T}$ is the portion of time that the system is on, and $T'$ satisfies the following:
\begin{align}
\begin{cases}
\eta_i(t_{k}+T')=e^{A_iT'}\eta_i(t_k)+A_i^{-1}(e^{A_iT'}-I)B_{on} \\
\eta_i^{(1)}(t_{k}+T')=u_i(t_k)-\delta/2 \nonumber \\
 \eta_i^{(1)}(t_k)=T_c^i(t_k). 
\end{cases}
\end{align}
$T_f^i(t_k)$ is the room temperature at $t_{k}+T$ given that the system is on during the entire period between $t_{k}$ and $t_{k}+T$, which satisfies the following:
\begin{align}
\begin{cases}
\eta_i(t_{k}+T)=e^{A_iT}\eta_i(t_k)+A_i^{-1}(e^{A_iT}-I)B_{on} \\
\eta_i^{(1)}(t_{k}+T)=T_f^i(t_k)\\
 \eta_i^{(1)}(t_k)=T_c^i(t_k).
\end{cases}
\label{eq:condition}
\end{align}
$T_f^i$ is defined in (\ref{eq:condition}) to characterize the condition in which the load is ``on'' for the entire period and therefore consumes the maximum energy. Intuitively, if the room temperature at $t_k$ is less than the lower bound of the control deadband ($T_c^i(t_k)\leq u_i(t_k)-\delta/2$), the power state will be ``off'' until the room temperature hits the boundary of the deadband. On the other hand, if $u_i(t_k)\leq T_f^i(t_k)+\delta/2$, it indicates that the load is always ``on'', and the room temperature does not hit the boundary for the entire period.

\begin{figure}[t]
\centering
\includegraphics[width = 0.7\linewidth]{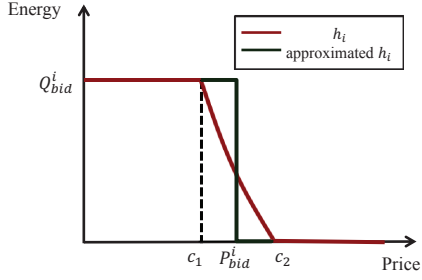}
\caption{The energy response curve $h_i$ and its approximation.}
\label{biddingana}
\end{figure}

To simplify the message space, we approximate $h_i$ with a step function as illustrated in Fig. \ref{biddingana}, where $c_1$ and $c_2$ are computed based on the control setpoint and user {\em type}. 
For notation convenience,  define $c_1=e_i(u_1,\theta_i)$ and $c_2=e_i(u_2,\theta_i)$ , where $u_1$ and $u_2$ are the temperature control setpoints corresponding to $c_1$ and $c_2$, respectively. For example, using the second-order ETP model (\ref{eq:inddynamics}) and control policy (\ref{eq:hvaccontrol}), $u_1$ and $u_2$ for the $i$th device can be obtained as:
\begin{align}
\begin{cases}
u_1&=T_c^i(t_k)+\delta/2 \\
u_2&=L A_i^{-1}e^{A_iT}(A_i\eta_i(t_k)+B^i_{on})-L A_i^{-1}B^i_{on}+ \delta/2 \\
  &=T_f^i(t_k)+\delta/2,
\end{cases}
\label{eq:bidding}
\end{align}
where $L=[1,0]$, and the power state of the $i$th TCL is ``on" at $t_k^-$. 

The step function in Fig. \ref{biddingana} can be fully characterized by two scalars: $P_{bid}^i$ and $Q_{bid}^i$, where $P_{bid}^i$ is the middle point of $c_1$ and $c_2$, while $Q_{bid}^i$ is the power consumption when the device is on during the market period. 
In this case, the message space of each user $M_i$ is reduced from a function space to a space of $\mathbb{R}^2_+$, and each bid is of the form $[P_{bid}^i,Q_{bid}^i]$.

\begin{remark}
Bidding and pricing can be viewed as information exchange between the coordinator and the loads that is essential for optimal decision making. Many advanced DLC methods also have communication requirements \cite{liu2013planning,han2010development,hao2013aggregate,TPS_13,lu2012evaluation,hao2013ancillary,burke2008robust} and can also accomplish certain group objectives. Some DLC strategies may even learn the user responses through the input/output user behaviors. The main difference of the proposed market-based approach lies in its emphasis on the quantitative incorporation of user preferences, the economic interpretation of user bids and coordination signals, and the encoding of internal load dynamics and user preference information into the bids.  
\end{remark}

\section{Output Based Bidding}
The proposed bidding strategy in Section IV assumes the knowledge of ETP model parameters. In practice these parameters are difficult to derive, and the ETP model used in the framework may be inaccurate in terms of characterizing the real energy consumption of the TCLs. To address these challenges, we present a joint state and parameter estimation framework, which enables users to compute bidding prices only based on local measurements.

In the ETP model (\ref{eq:inddynamics}), $A_i$ is a constant, while $B_{on}^i$ and $B_{off}^i$ are time varying, which depend on the outside temperature and the solar heat gain. Let $\zeta_k\in R^2$ be a vector denoting the outside temperature and the solar heat gain. We assume that $\zeta_k$ can be measured or estimated. When we have some rough statistical information about the model parameters, the system dynamics can be captured with an uncertain discrete dynamic model with Gaussian noise::
\begin{align}
\begin{cases}
&\eta_i(t_{k}) =
\bar{A}_i\eta_i(t_{k-1})+\bar{B}_i \zeta_{k-1}+ \bar{C}_i+w_{k-1}^i \\
&y_i(t_k) =L\eta_i(t_k)+v_{k}^i  ,
\label{eq:uncertainmodel}
\end{cases}
\end{align}
where $y_i(t_k)$ is the output measurement (air temperature), $L=[1,0]$ and we have two linear subsystems depending on the power state of the load:
\begin{align}
\bar{C}_i =\begin{cases}
\bar{C}_{on}^i  & \text{ if } q_i(t_{k-1})=1\\
\bar{C}_{off}^i  & \text{ if } q_i(t_{k-1})=0.
\end{cases} 
\end{align}
In this model, the dependence of the load dynamics on the external signal $\zeta_k$ is made explicit. Therefore, $A_i$, $\bar{B}_i$, $\bar{C}_{on}^i$ and $\bar{C}_{on}^i$ are time invariant unknown parameters. Here we assume that all the noise terms follow the Gaussian distributions:
\begin{equation}
\label{eq:uncertainties}
\begin{cases}
w_k^i\sim \mathcal{N}(w_k^i\mid 0,\Omega_i)   \\
v_k^i \sim \mathcal{N}(v_k^i\mid 0,\Sigma_i) \\
\mu_i \sim \mathcal{N}(\mu_i\mid 0,\Phi_0^i)  .
\end{cases}
\end{equation}
Let $\eta_i(t_1)=m_0^i+\mu_i$ be the initial state ($\mu_i$ is the noise). Denote $\sigma_i=[\bar{A}_i,\bar{B}_i,\bar{C}_{on}^i,\bar{C}_{off}^i,\Omega_i,\Sigma_i, m_0^i, \Phi_0^i]$ as the unknown parameter to be estimated. The problem can be then formulated as estimating $\sigma_i$ base on local measurements, $Y_i=(T_c^i(t_1),\ldots,T_c^i(t_M))$. This can be cast as a joint state-parameter estimation problem, which can be solved using the  expectation maximization (EM) algorithm \cite[chap. 13]{bishop2006pattern}.
The EM algorithm is a two-stage iterative optimization technique for finding the maximum likelihood solution for the unknown parameters. In other words, it finds the optimal $\sigma_i$ that maximizes the likelihood function $p(Y_i|\sigma_i)$, where $Y_i=(y_i(t_1),\ldots,y_i(t_M))$. The EM algorithm starts from some initial selection for the model parameters, $\sigma_{old}$. In the first stage (E-step), we evaluate the posterior distribution of the state $p(Z_i|Y_i,\sigma_{old})$ assuming that all the parameters are known, where $Z_i=(\eta_i(t_1),\ldots,\eta_i(t_M))$. In the second stage (M-step), the derived  posterior distribution is used to find the updates of $\sigma_i$ that maximizes the expectation of the logarithm of the complete-data likelihood function, which is:
\begin{equation}
Q(\sigma_i,\sigma_{old})=\sum_{Z_i} p(Z_i|Y_i,\sigma_{old})\text{ln}\left(p(Y_i,Z_i|\sigma_i)\right) .
\label{completelikelihood}
\end{equation}
After the update for the parameter estimation is derived in the M step, we assign it to $\sigma_{old}$ and go back to E-step. This procedure is iterated until the estimation of the state and parameters converges. 
The detailed steps of the proposed estimation algorithm for the ETP model can be briefly summarized as follows:

\subsubsection{The E Step}
The E step finds the distribution for the system state $\eta_i(t_k)$ conditioned on the full observation sequence, $Y_i=(y_i(t_1),\ldots,y_i(t_M))$, assuming that the model parameters are known as $\sigma_{old}$. This inference problem can be solved efficiently using the sum-product algorithm \cite{bishop2006pattern} in two steps: first, the distribution of state $\eta_i(t_k)$ conditioned on a partial observation sequence $(y_i(t_1),\ldots,y_i(t_k))$ can be derived with a Kalman filter; second, the conditional distribution $p(\eta_i(t_k)\mid Y_i)$ can be found with a Kalman smoother.

Denote $\hat{\alpha}(\eta_i(t_k))$ as the conditional distribution $p(\eta_i(t_k)\mid y_i(t_1),\ldots,y_i(t_k))$, which satisfies:
 \begin{equation}
\label{eq:14}
\hat{\alpha}(\eta_i(t_k))=\mathcal{N}(\eta_i(t_k)\mid \mu_k, \Phi_k)  ,
\end{equation}
where $\mathcal{N}$ stands for Gaussian distribution with $u_k$ and $\Phi_k$ as its mean and covariance, respectively. 
In the context of linear-Gaussian systems, the sum-product algorithm \cite{bishop2006pattern} gives the following recursion equations:

\begin{align}
\begin{cases}
\mu_k=&\bar{A}_i\mu_{k-1}+\bar{B}_i \zeta_{k-1}+\bar{C}_i+K_k\Big(y_i(t_k)-L\bar{A}_i\mu_{k-1} \\
&-L\bar{B}_i\zeta_{k-1}-L\bar{C}_i\Big) \\
\Phi_k=&(I-K_kL)P_{k-1} , \\
\end{cases}
\end{align}
where $\bar{C}_i=\bar{C}_{\text{on}}^i$ if the $i$th load is ``on", and $\bar{C}_i=\bar{C}_{\text{off}}^i$ if otherwise.  $P_k$ and the Kalman gain matrix is defined as:
 \begin{equation}
\label{eq:18}
\begin{cases}
P_{k-1}=\bar{A}_i \Phi_{k-1}\bar{A}_i^T+\Omega_i \\
K_k=P_{k-1}L^T (LP_{k-1}L^T+\Sigma_i)^{-1} .
\end{cases}
\end{equation}
The initial conditions for the recursion equation are given by:
\begin{equation}
\label{eq:20}
\begin{cases}
\mu_1=m_0+K_1(y_i(t_1)-Lm_0) \\
\Phi_1=(I-K_1L)\Phi_0^i , 
\end{cases}
\end{equation}
where $K _1=\Phi_0^iL^T(L\Phi_0^iL^T+\Sigma_i)^{-1}$. 
 
With the above recursion equations,  we can derive the distribution for $\eta_i(t_k)$ conditioned on the observations from $y_i(t_1)$ to $y_i(t_k)$. Next we turn to the problem of finding the probability distribution for $\eta_i(t_k)$ given all observations from $y_i(t_1)$ to $y_i(t_M)$. Denote the conditional distribution  $p(\eta_i(t_k)\mid Y_i)$ as $\gamma(\eta_i(t_k))$, which satisfies:
 \begin{equation}
\label{eq:21}
\gamma(\eta_i(t_k))=\mathcal{N}(\eta_i(t_k)\mid \hat{\mu}_k, \hat{\Phi}_k) .
\end{equation}
The sum-product algorithm gives the following recursion equations:
\begin{equation}
\label{eq:23}
\begin{cases}
\hat{\mu}_k=\mu_k+J_k(\hat{\mu}_{k+1}-\bar{A}_i\mu_k-\bar{B}_i\zeta_{k}-\bar{C}_i) \\
\hat{\Phi}_k=\Phi_k+J_k(\hat{\Phi}_{k+1}-P_k)J_k^T ,
\end{cases}
\end{equation}
where $J_k=\Phi_k\bar{A}_i^T(P_k)^{-1}$. 

With the recursion equation presented above, the conditional distribution $p(\eta_i(t_k)\mid Y_i)$ can be computed using backward induction.

\subsubsection{The M Step}
The M step tries to find the parameter update that maximizes the logarithm of the complete-data likelihood function (\ref{completelikelihood}). Equation (\ref{completelikelihood}) indicates that aside from the conditional distribution $p(\eta_i(t_k)\mid Y_i, \theta_{old})$ (already obtained in the E step), the likelihood function also depends on the joint distribution $p(Z_i,Y_i\mid \sigma_i)$.
For the linear-Gaussian system (\ref{eq:uncertainmodel}), the logarithm of this joint distribution $p(Z_i,Y_i\mid \sigma_i)$ is given by:
\begin{align}
\text{ln}\,p(Z_i,Y_i\mid \sigma_i) = &\sum_{k=2}^M \text{ln}\,p(\eta_i(t_k)\mid \eta_i(t_{k-1}),\bar{A}_i,\bar{B}_i,\bar{C}_i,\Omega_i) \nonumber \\
&+\sum_{k=1}^M \text{ln}\,p(y_i(t_k)\mid \eta_i(t_k),L,\Sigma_i) \nonumber \\
&+\text{ln}\,p(\eta_i(t_1)\mid m_0^i,\Phi_0^i) ,  \label{eq:7}
\end{align}
where the dependence of the joint distribution on the unknown model parameters is  made explicit. 
The complete-data likelihood function $Q(\sigma_i,\sigma_{old})$ can be then obtained by taking the expectation of (\ref{eq:7}) over $Z_i$ using the posterior distribution $p(\eta_i(t_k)\mid Y_i, \theta_{old})$ derived in the E step. 

Let $\sigma_i'=(\bar{A}_i',\bar{B}_i',\bar{C}_{i,on}',\bar{C}_{i,off}',\Omega_i',\Sigma_i', m_0',\Phi_0')$ be the update of the unknown parameter in the M step, i.e., $\sigma_i'=\argmax_{\sigma_i} Q(\sigma_i,\sigma_{old})$. 
The explicit formula for each component of $\sigma_i'$ is given as follows (please refer to \cite{bishop2006pattern} for the detailed derivation):
\begin{enumerate}
\item Maximizing (\ref{completelikelihood}) over $m_0^i$ and $\Phi_0$, the updates can be derived as: 
\begin{equation*}
\begin{cases}
m_0'=\mathbb{E}[\eta_i(t_1)] \\
\Phi_0'=\mathbb{E}[\eta_i(t_1)\eta_i(t_1)^T]-\mathbb{E}[\eta_i(t_1)]\mathbb{E}[\eta_i(t_1)^T]  .
\end{cases}
\end{equation*}

\item Maximizing the likelihood function (\ref{completelikelihood}) over $\bar{A}_i$, the update of $\bar{A}_i$ is given by:
\begin{align}
\hspace{-0.65cm}
\bar{A}_i'&=\Big( \sum_{k=2}^M \mathbb{E}[\eta_i(t_k)\eta_i(t_{k-1})^T]-\bar{B}_i\zeta_{k-1}\mathbb{E}[\eta_i(t_{k-1})^T] \nonumber \\ 
&-\bar{C}_i \mathbb{E}[\eta_i(t_{k-1})^T] \Big) \times	
 \Big( \sum_{k=2}^M \mathbb{E}[\eta_i(t_{k-1})\eta_i(t_{k-1})^T]\Big)^{-1}. \label{updateA}
\end{align}

\item Maximizing the likelihood function (\ref{completelikelihood}) over $\bar{B}_i$, we can derive the update of $\bar{B}_i$ as follows:
\begin{align}
\label{updateB}
\hspace{-0.65cm}
\bar{B}_i'=&\left(\sum_{k=2}^M \Big( \mathbb{E}[\eta_i(t_k)]-\bar{A}_i'\mathbb{E}[\eta_i(t_{k-1})]-\bar{C}_i\Big)\zeta_{k-1}^T  \right)\times \nonumber \\
&\left( \sum_{k=2}^M \zeta_{k-1}\zeta_{k-1}^T \right)^{-1}
\end{align}

\item Maximizing the likelihood function (\ref{completelikelihood}) over $\bar{C}_{on}^i$ and $\bar{C}_{off}^i$, the updates are given as:
\begin{align}
\label{eq:11}
\hspace{-0.65cm}
\left\{
\begin{aligned}
&\bar{C}_{i,\text{on}}'=\dfrac{1}{\vartheta_{1}}\underset{{\scriptstyle k\in M_{1}}}{\sum}\Big(\mathbb{E}[\eta_i(t_k)]-\bar{A}_i'\mathbb{E}[\eta_i(t_{k-1})]-\bar{B}_i'\zeta_{k-1}\Big) \\
&\bar{C}_{i,\text{off}}'=\dfrac{1}{\vartheta_{2}}\underset{{\scriptstyle k\in M_{2}}}{\sum}\Big(\mathbb{E}[\eta_i(t_k)]-\bar{A}_i'\mathbb{E}[\eta_i(t_{k-1})]-\bar{B}_i'\zeta_{k-1}\Big),
\end{aligned}
\right.
\end{align}
where $M_{1} \subseteq \{1,2,\ldots,M\}$ denotes the time instants when the system is on, and $\vartheta_{on}$ represents the size of $M_{1}$. $M_2$ and $\vartheta_{2}$ are defined similarly.

\item The update function for $\Omega_i$ can also be derived by maximizing the likelihood function with respect to $\Omega_i$, which gives:

\begin{align}
\hspace{-1.3cm}
\Omega_i'=&\frac{1}{M-1}\sum_{k=2}^M
\Big\{\bar{A}_i'\mathbb{E}[\eta_i(t_{k-1})\eta_i(t_{k-1})^T]\bar{A}_i'^T  \nonumber \\
&-\mathbb{E}[\eta_i(t_k)\eta_i(t_{k-1})^T]\bar{A}_i'-\bar{A}_i'\mathbb{E}[\eta_i(t_{k-1})\eta_i(t_k)^T] \nonumber \\
&+\mathbb{E}[\eta_i(t_k)\eta_i(t_k)^T]-(\bar{B}_i'\zeta_{k-1}+\bar{C_i}')\mathbb{E}[\eta_i(t_{k})^T] \nonumber \\
&-\mathbb{E}[\eta_i(t_{k})](\zeta_{k-1}^T\bar{B}_i'^T+\bar{C}_i'^T)+\bar{C_i}'\mathbb{E}[\eta_i(t_{k-1})^T]\bar{A}_i'^T \nonumber \\
&+\bar{B}_i'\zeta_{k-1}\mathbb{E}[\eta_i(t_{k-1})^T]\bar{A}_i'^T+\bar{A}_i'\mathbb{E}[\eta_i(t_{k-1})]\zeta_{k-1}^T\bar{B}_i'^T \nonumber \\
&+ \bar{A}_i'\mathbb{E}[\eta_i(t_{k-1})]\bar{C}_i'^T+(\bar{B}_i'\zeta_{k-1}+\bar{C}_i')\zeta_{k-1}^T\bar{B}_i'^T \nonumber \\
&+ (\bar{B}_i'\zeta_{k-1}+\bar{C}_i')\bar{C}_i'^T  \Big\}   . \label{eq:updategamma}
\end{align}

\item  The update for $\Sigma_i$ can also be obtained similarly:

\begin{align}
\Sigma_i'=&\frac{1}{M}\sum_{k=1}^M
\{y_i(t_k)y_i(t_k)^T-L\mathbb{E}[\eta_i(t_k)]y_i(t_k)^T- \nonumber \\
&y_i(t_k)\mathbb{E}[\eta_i(t_k)^T]L+L\mathbb{E}[\eta_i(t_k)\eta_i(t_k)^T]L \} . \label{eq:13}
\end{align}

\end{enumerate}

In the above update equations, $\mathbb{E}[\eta_i(t_k)]$ and $\mathbb{E}[\eta_i(t_k)\eta_i(t_{k})^T]$ can be computed based on the conditional distribution $p(\eta_i(t_k)\mid Y_i)$ obtained in the E step, while the pairwise expectation $\mathbb{E}[\eta_i(t_k)\eta_i(t_{k-1})^T]$ can be derived using Bayesian Theorem. The expressions for these expectations are as follows \cite{bishop2006pattern}: 
\begin{equation}
\label{eq:5}
\begin{cases}
\mathbb{E}[\eta_i(t_k)]&=\hat{\mu}_{k} \\
\mathbb{E}[\eta_i(t_k)\eta_i(t_{k-1})^T]&=J_{k-1}\hat{\Phi}_k+\hat{\mu}_{k}\hat{\mu}_{k-1} \\
\mathbb{E}[\eta_i(t_k)\eta_i(t_{k})^T]&=\hat{\Phi}_k+\hat{\mu}_k \hat{\mu}^T_k .
\end{cases}
\end{equation}
After the update $\sigma_i'$ of the the estimated parameter is derived, we assign it to $\sigma_{old}$ and go back to the E step. This process is repeated until the update for the estimated parameters converges. When the estimation for $\sigma_i$ is obtained, each user can compute the bidding prices based on (\ref{eq:bidding}) and the bidding curve as shown in Fig. \ref{fig:bidding}. The output-based bidding algorithm is summarized as Algorithm 1.

\begin{algorithm}[tph]
\caption{The output-based bidding algorithm} \label{alg:1}
\begin{algorithmic}[1]

\REQUIRE Initial guess of parameters $\sigma_{old}$ and measurement sequence $Y_i$.

\WHILE {convergence criteria not satisfied}
\STATE  Find the conditional distribution of system state $p(Z_i|Y_i,\sigma_{old})$ using Kalman Filter.

\STATE Find $\sigma_i$ that maximizes the complete-data likelihood function $Q(\sigma_i,\sigma_{old})=\sum_{Z_i} p(Z_i|Y_i,\sigma_{old})\text{ln}p(Y_i,Z_i|\sigma_i)$. Denote the optimal solution as $\sigma_{new}$

\STATE Update the parameter estimation: $\sigma_{old}=\sigma_{new}$.

\ENDWHILE

\ENSURE The estimation of the state trajectory $Z_i$ and the parameters $\sigma_i$.

\end{algorithmic}
\end{algorithm}

\section{Case Studies}
This section applies the proposed market mechanism and the learning scheme to the TCL coordination problem considered in the GridWise$\textsuperscript{\textregistered}$ demonstration project \cite{fuller2011analysis}, and presents simulation results to demonstrate the effectiveness of the proposed market mechanism.

\subsection{Simulation Setup}
Similar to the GridWise$\textsuperscript{\textregistered}$ project, we consider a realistic scenario where  each user is equipped with a smart thermostat that can measure the room temperature and communicate to the coordinator. At each period, the device measures the current room temperature and submits a bidding price based on a bidding curve. The coordinator collects all the user bids and clears the energy market with a price. Each device will then determine the control setpoint in response to this energy price, which modifies the load dynamics and affects the bids for the next period. This framework is validated in Matlab using the parameters generated in GridLab-D \cite{chassin2008gridlab}. The details of the demonstration project can be found in Section II.

The second-order ETP model is used to capture the load dynamics of the TCLs. The ETP model parameters depend on various building parameters, such as glass type, floor area, area per floor,  glazing layers and material,  etc. Detailed description of these parameters and their relations to the ETP model parameters can be found in 
\cite{gridlab}. In the simulation, 1000 sets of building parameters are generated. A few important parameters are randomly generated using the same approach as in \cite{TPS_13}, while the rest take their default values from GridLAB-D. 
Throughout the simulation, the aggregated power of the unresponsive loads is assumed to be $12MW$, and the feeder power capacity is $15MW$. In addition, we refer to the ``power" consumption of the load as the average power during one market period (5 minute), unless otherwise stated.

We use the weather data and the Typical Meteorological Year (TMY2) data for Columbus, OH, obtained from \cite{weatherdata}, \cite{marion1995user}, which includes air temperature and solar radiation. The wholesale energy price is from the PJM market \cite{PJMLMP}. It is modified to a retail rate in $\$/kWh$ plus a retail modifier as defined  by AEP's tariff \cite{aeptariff}, and we define this retail price as the \emph{base price}.

\begin{figure}[t]
\centering
\includegraphics[width = 0.80\linewidth]{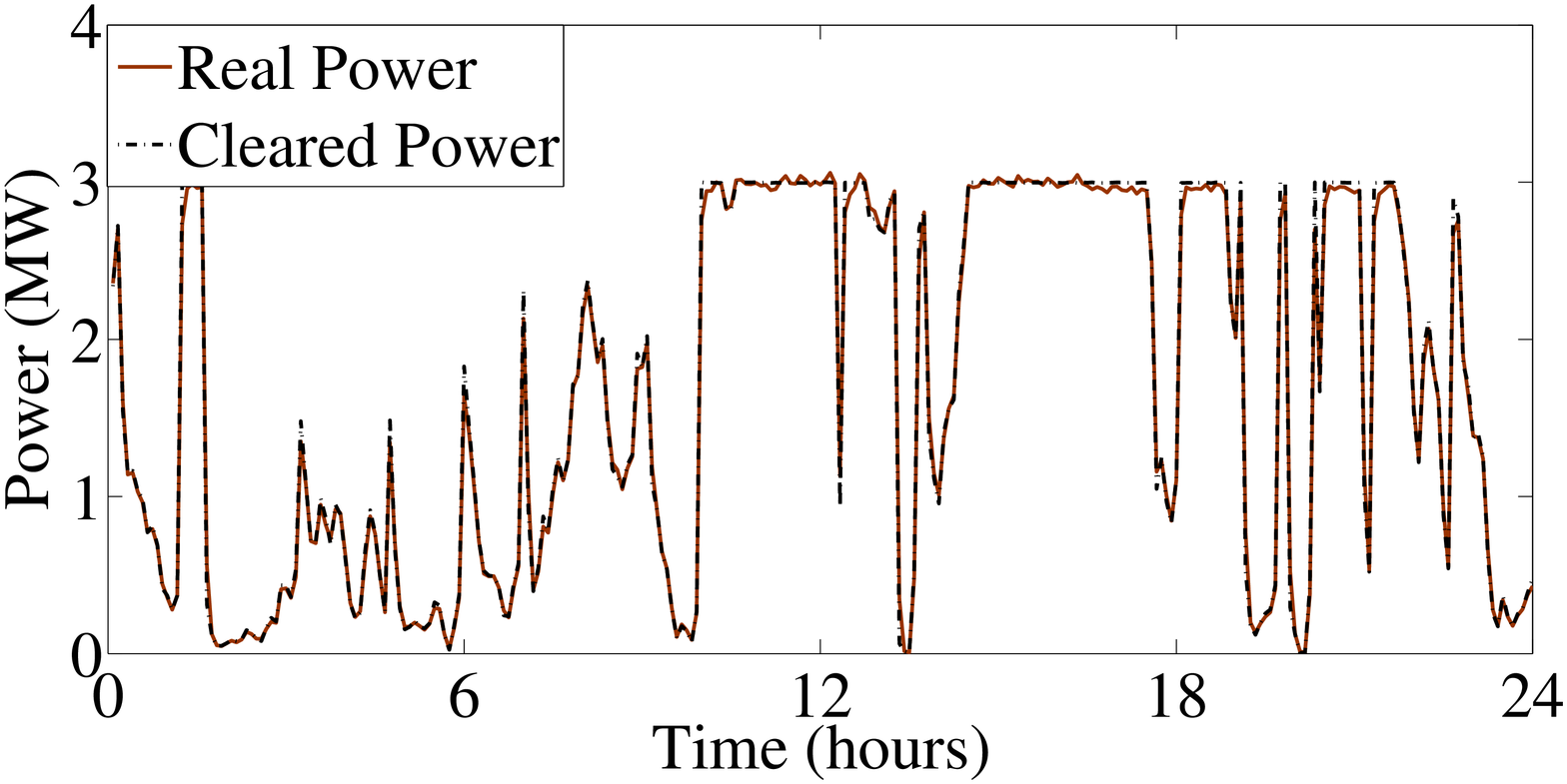}
\caption{Comparison of the actual power trajectory and the cleared power. The outside air temperature record is on August 20, 2009 in Columbus, OH.}
\label{powertrajectory1}
\end{figure}

\begin{figure}[t]
\centering
\includegraphics[width = 0.85\linewidth]{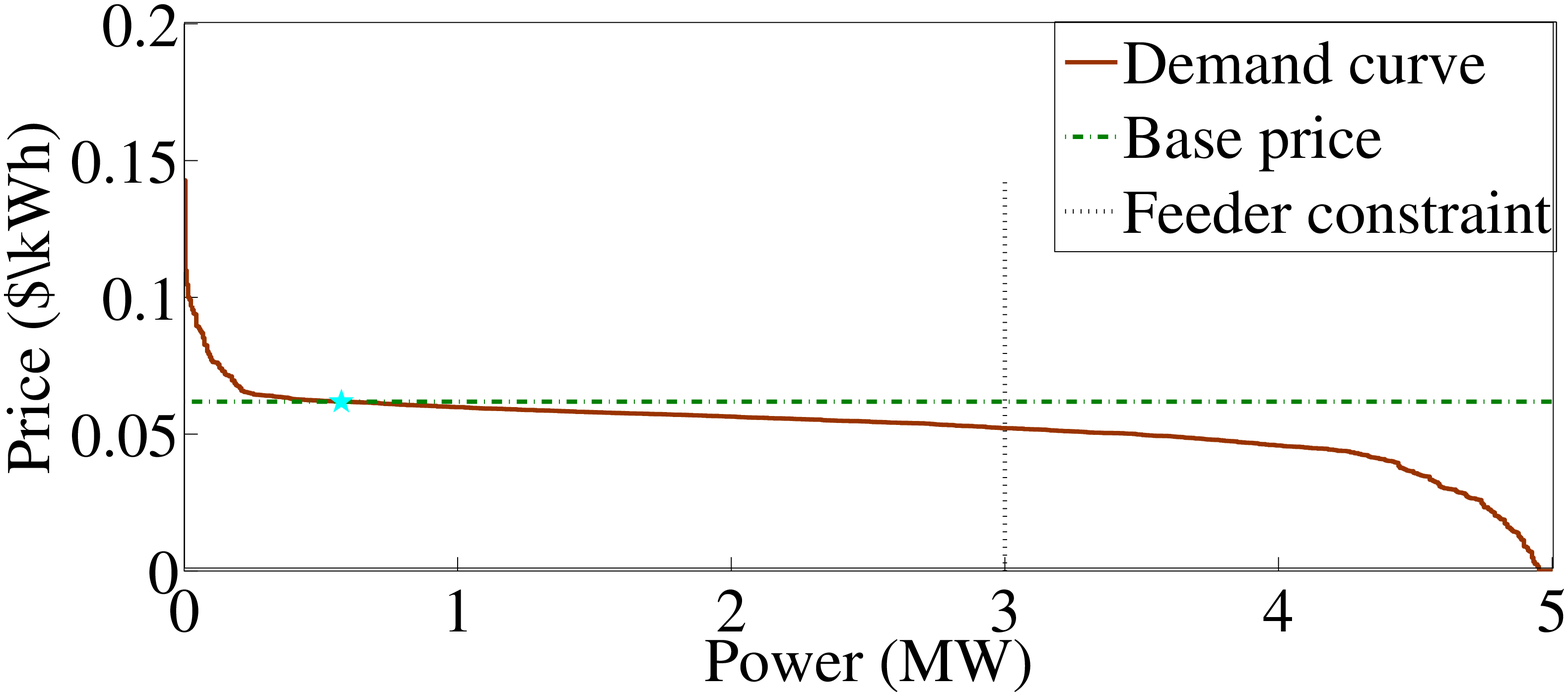}
\caption{The demand curve and the market clearing process at 08:20 AM. When the total demand is less than the feeder capacity constraint, the market price is equal to the base price.}
\label{clearnocon}
\end{figure}

\begin{figure}[t]
\centering
\includegraphics[width = 0.85\linewidth]{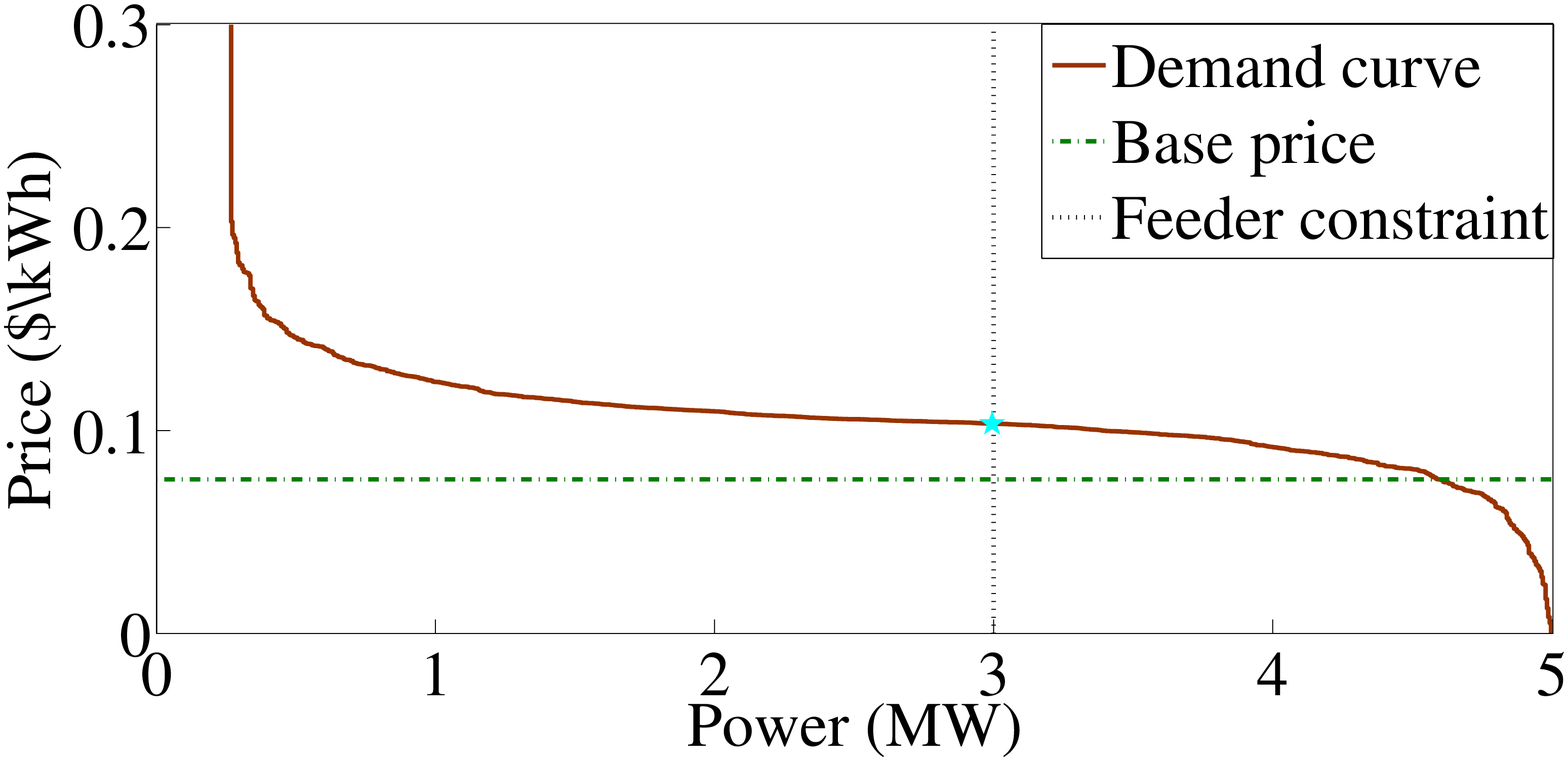}
\caption{The demand curve and the market clearing process at 04:40 PM. When the total demand exceeds the feeder capacity constraint, the market price is higher than base price to respect the feeder capacity constraint.}
\label{clearcon}
\end{figure}

\begin{figure}[t]
\centering
\includegraphics[width = 0.85\linewidth]{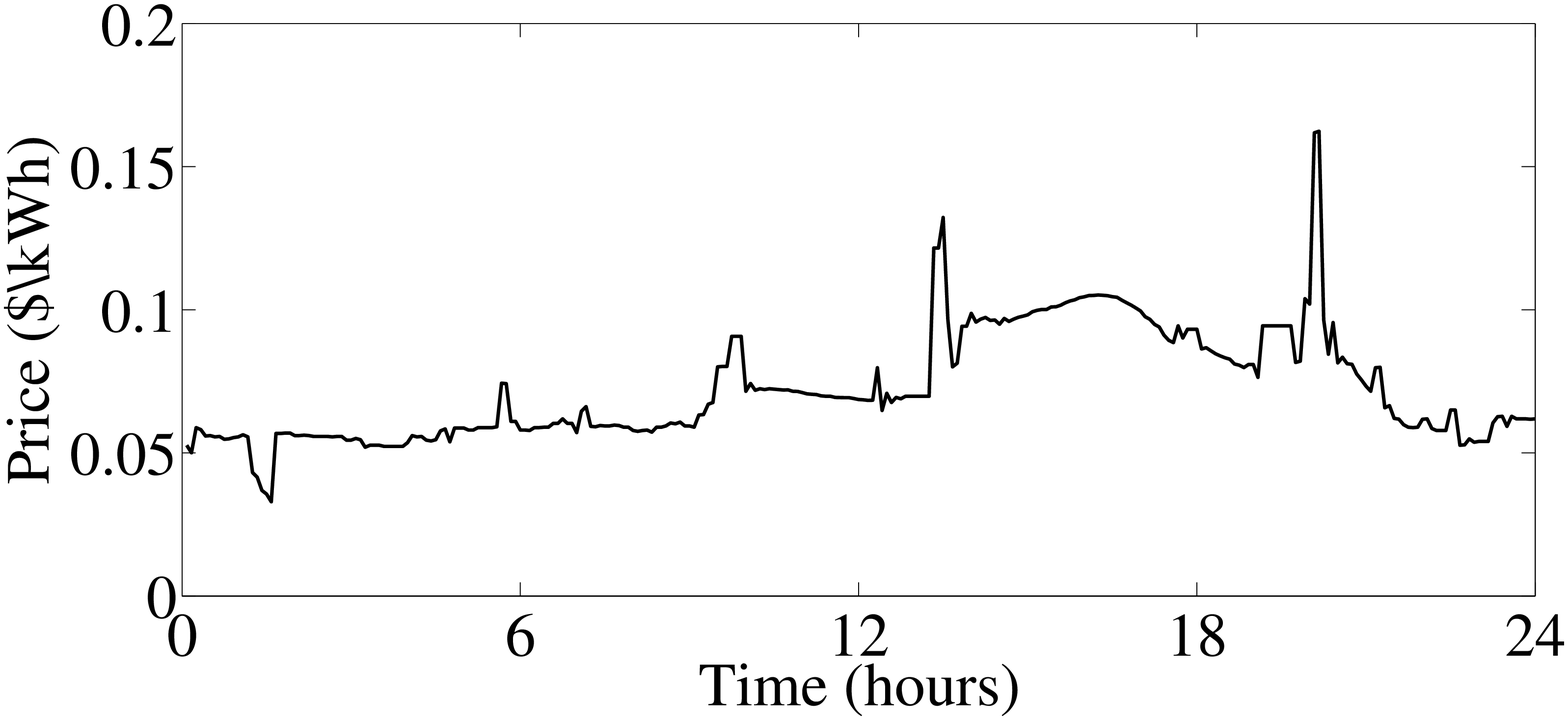}
\caption{Comparison of the trajectories of the market clearing price and wholesale energy price. The price is higher during the congestion, which effectively caps the peak energy at key times.}
\label{pricecurve1}
\end{figure}

\subsection{Results and Analysis}
First the proposed mechanism is evaluated in the deterministic case, where each user can accurately estimate the unknown parameters. Each user submits a bid as described in Fig. \ref{fig:bidding}, and the market is cleared according to the proposed outcome function (\ref{EQ_1})-(\ref{EQ_3}). Air temperature record is on August 20 (mild day), 2009 in Columbus, OH, and the aggregated power trajectory is presented in Fig. \ref{powertrajectory1}.

It shows that the power trajectory is effectively capped below the feeder power capacity for the entire day. Notice that whenever the coordinator clears the market with an energy price, there is a corresponding power on the aggregated demand curve (as shown in Fig. \ref{fig:demandcurve1} and Fig. \ref{fig:demandcurve2}). We call it the cleared power, which stands for the coordinator's estimation on the aggregated power consumption before the market is cleared. Simulation results demonstrate that this cleared power accurately matches the actual power consumption. This enables the coordinator to select proper prices to effectively achieve the desired aggregated power consumption. To demonstrate how it works, we randomly choose two market periods and present their market clearing procedures in Fig. \ref{clearnocon} and Fig. \ref{clearcon}, respectively. When there is no power congestion, the coordinator can directly pass the base price to individual users. This case is shown in Fig. \ref{clearnocon}, which corresponds to 08:20 AM in Fig. \ref{powertrajectory1}. When the power congestion occurs, the coordinator should clear the market at the intersection of the aggregated demand curve and the curve of the feeder power constraint. This case is presented in Fig. \ref{clearcon}, which corresponds to 04:40 PM. 

\begin{figure}[t]
\centering
\includegraphics[width = 0.85\linewidth]{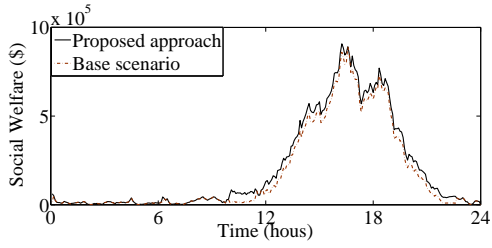}
\caption{Comparison of the social welfare of the proposed pricing strategy and the base scenario. The base scenario adopt RTP and multiplies the base price by a fixed ratio to cap the total energy.}
\label{compare}
\end{figure}

\begin{figure}[t]
\centering
\includegraphics[width = 0.75\linewidth]{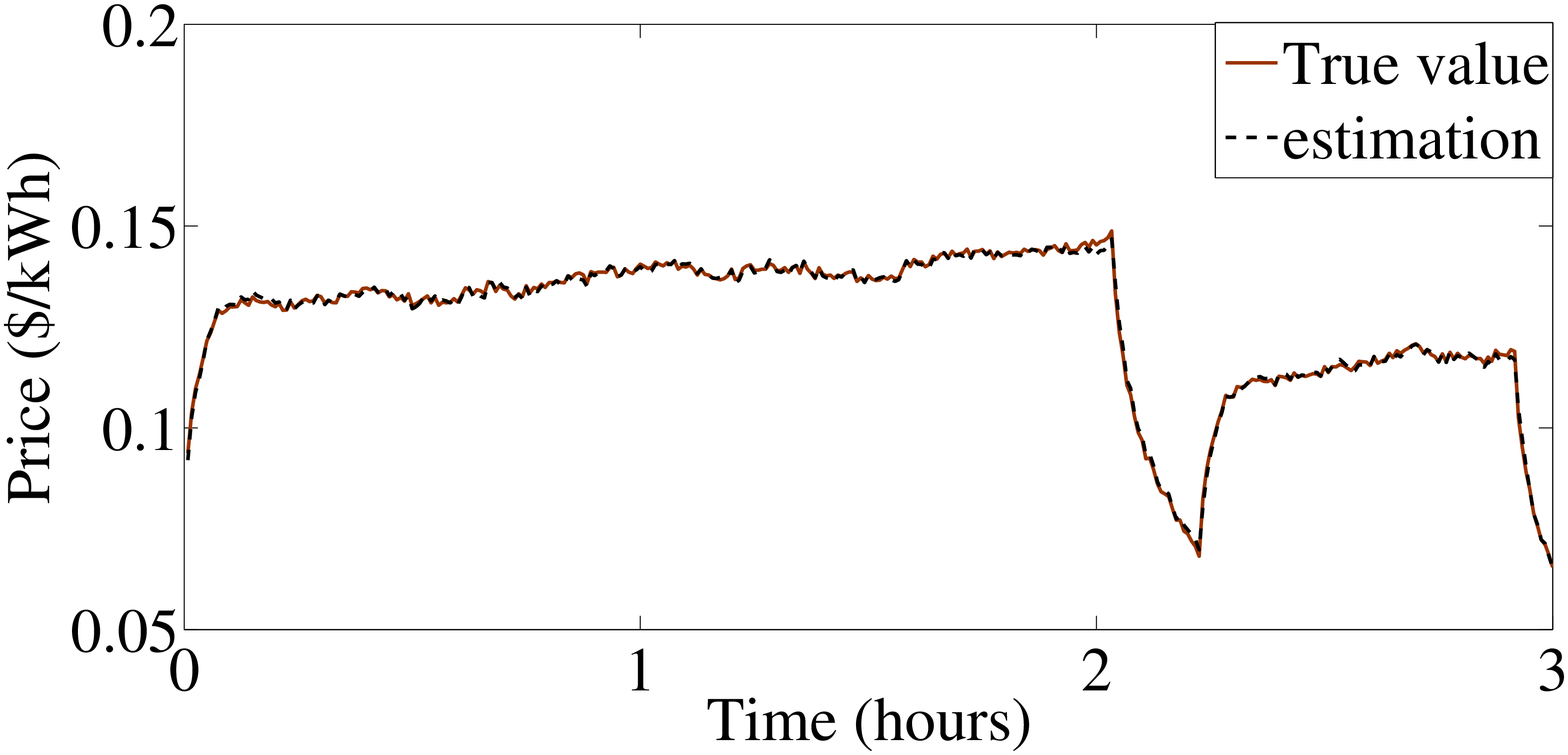}
\caption{The estimation result of the output-based bidding algorithm. The initial guess of the simulation is randomly selected from $90\%$ and $110\%$ of its true value.}
\label{EM}
\end{figure}

\begin{figure}[t]
\centering
\includegraphics[width = 0.85\linewidth]{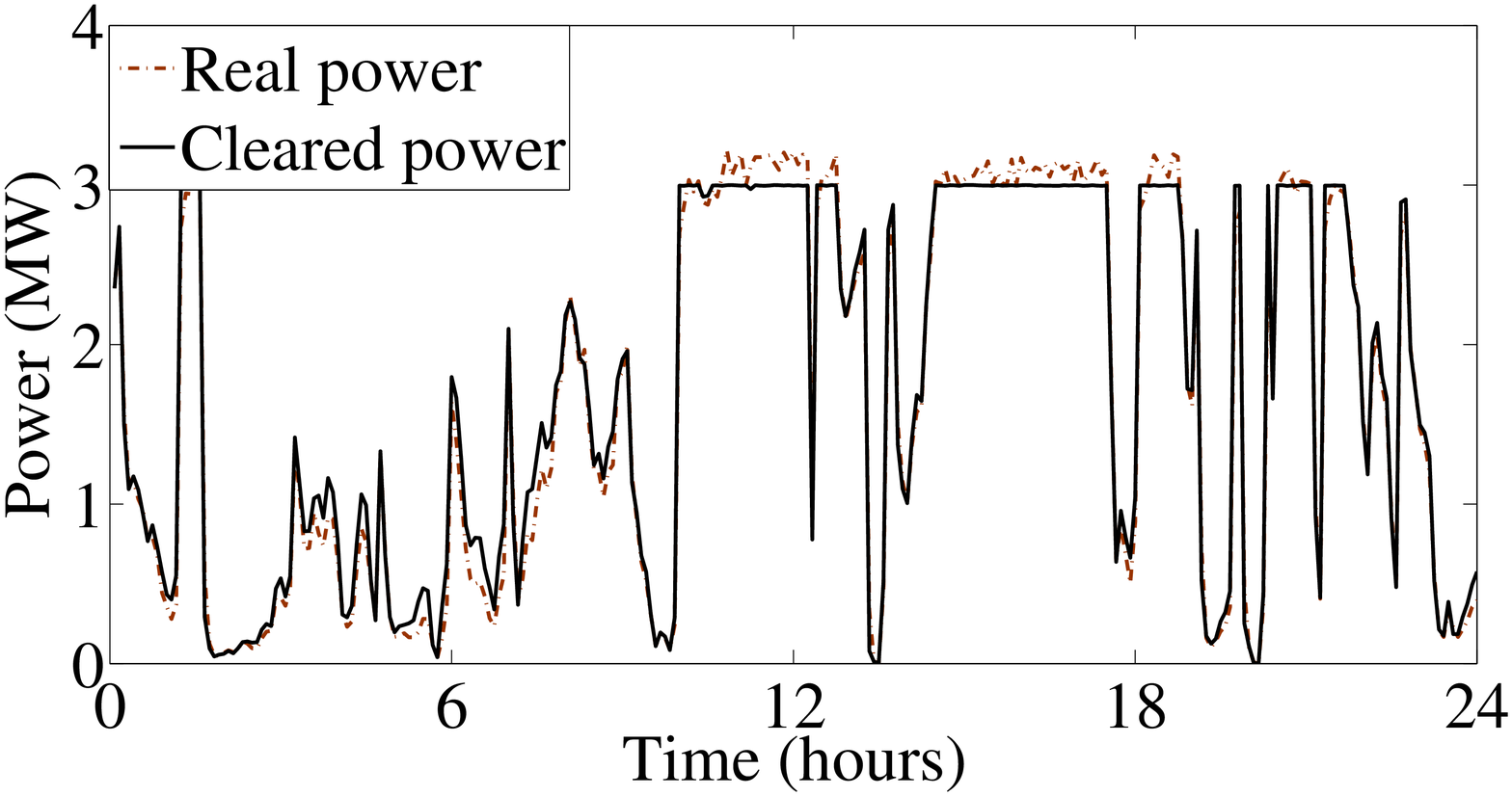}
\caption{The actual power trajectory and the cleared power based on the demand curve with $2\%$ bidding error. The outside air temperature record is on August 20, 2009 in Columbus, OH.}
\label{powertrajectory3}
\end{figure}

\begin{figure}[t]
\centering
\includegraphics[width = 0.85\linewidth]{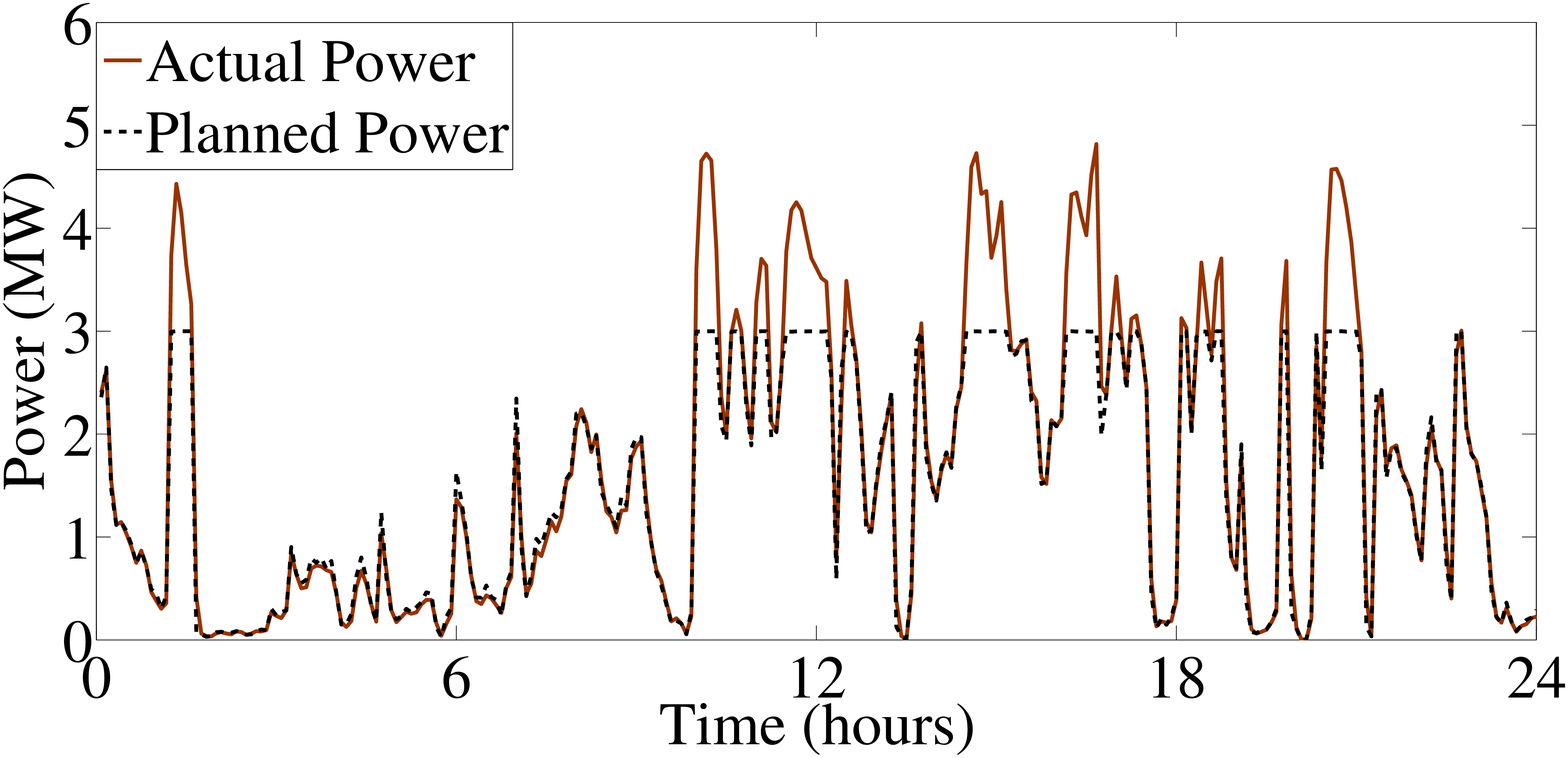}
\caption{The aggregated power and the cleared power under real time pricing strategy. The base price is directly passed to the retail market.}
\label{RTP}
\end{figure}

The  trajectories of the market clearing price and the base price are shown in Fig. \ref{pricecurve1}. 
The figure shows that the market clearing price is equal to the base price when there is no congestion, and is higher than the base price during congestion hours. Notice that when $C'(\alpha)\neq P_{base}$, the market price can be different from base price in uncongested period as well, where $\alpha$ is the total energy purchased from the wholesale market. In addition, a few price spikes can be observed in the simulation result. However, these spikes are not created by the proposed framework, but instead caused by the  fluctuations of the base price. 

Furthermore, to evaluate the proposed mechanism in terms of social welfare maximization, we compare it with a base scenario, where Real Time Pricing (RTP) is adopted to cap the power  in a heuristic way. More specifically, when there is no congestion, the market clearing price is equal to the base price. When the power congestion occurs, the clearing price is the base price multiplied by a fixed ratio $\gamma$, which is greater than 1 and can cap the aggregated power below the limit effectively. We can run simulations to find such ratios, and among all the possible ratios, we choose the minimum one that can cap the aggregated power below the feeder capacity. Since the social welfare of the two scenarios will be the same during the uncongestion period ($\gamma=1$), we use the weather data on August 16, 2009, where more power congestion can be observed due to the elevated temperature. In this case $\gamma=2.6$, and the social welfare of the two pricing strategies is shown in Fig. \ref{compare}. The simulation results demonstrate that the proposed optimal pricing strategy always outperforms the base scenario in terms of social welfare. Notice that in this paper we run simulations for different values of $\gamma$
 to find the minimum value that can cap the aggregated power. However, this ratio is difficult to derive in practice, and therefore a much more conservative value has to be used to operate the power grid safely. This will further reduce the social welfare of the base scenario.

\subsection{The Output-based Bidding Algorithm}
\begin{figure}[t]
\centering
\includegraphics[width = 0.85\linewidth]{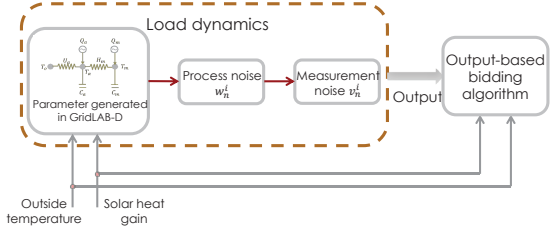}
\caption{Illustration for the out-based bidding algorithm. The ETP model parameters are default values in GridLAB-D generated according to a non-Gaussian distribution.}
\label{diagram}
\end{figure}

This subsection shows how the proposed output-based bidding algorithm can be used to accurately estimate the bidding prices. Fig. \ref{diagram} shows the simulation setups for the proposed algorithm. In the simulation, the ETP model parameters are the default values in GridLab-D, which are generated based on some non-Gaussian distribution, while the process noises  $w_n^i$ and the measurement noises $v_n^i$ are Gaussian. In addition, we assume that each device can locally measure its room temperature every minute, and store all the measurements for the past 6 hours, in which case $M=360$. The algorithm is started with an initial guess $\sigma_{old}$ with $10\%$ error. In other words, each element of the initial guess $\sigma_{old}$ is generated by randomly selecting a value between $90\%$ and $110\%$ of its true value. With the estimated parameters derived in the output-based bidding algorithm, each user can compute the bidding prices based on the realistic bidding strategies shown in Fig. \ref{biddingana}.  
Here we choose a random user in a random market period and present its estimation result in Fig. \ref{EM}. In this figure, the estimated bidding price is derived from the output-based bidding algorithm, while the true bidding price is computed based on the true value of the unknown parameters. It can be seen that the estimated bid closely follows the true value, with an average estimation error less than $1\%$. 

When all the users apply the output-based bidding algorithm to compute the biding prices, an error (of less than $1\%$) will be introduced. Now we evaluate how this estimation error can affect the aggregated power response.  To implement the estimation framework, each device will locally perform the output-based bidding algorithm during each market period, which just takes $5.5$ seconds on a laptop with 2.5GHz Intel i5 processor and 8G memory. However, it is computationally intensive to do the centralized simulations for all the users over 24 hours to show how the estimation error affects the aggregated power response. For this reason, instead of directly incorporating the output-based bidding algorithm in individual simulations, we add a simulated error of $2\%$ (this simulated error is larger than the actual error of the output-based bidding algorithm) to each user's bidding price. The simulation results with this bidding error are presented in Fig. \ref{powertrajectory3}. It can be seen that the aggregated power is effectively capped below the feeder capacity during $84\%$ of the time. In the cases where the feeder capacity constraint is violated, the aggregated power exceeds the power limit by $1.1\%$ on average, and the maximum violation occurred at 4:15 PM, where the power limit is exceeded by $3.3\%$.  Notice that these small violations can be easily fixed by adjusting the feeder capacity constraint in the formulation to be slightly more conservative.


\begin{figure}[t]
\centering
\includegraphics[width = 0.85\linewidth]{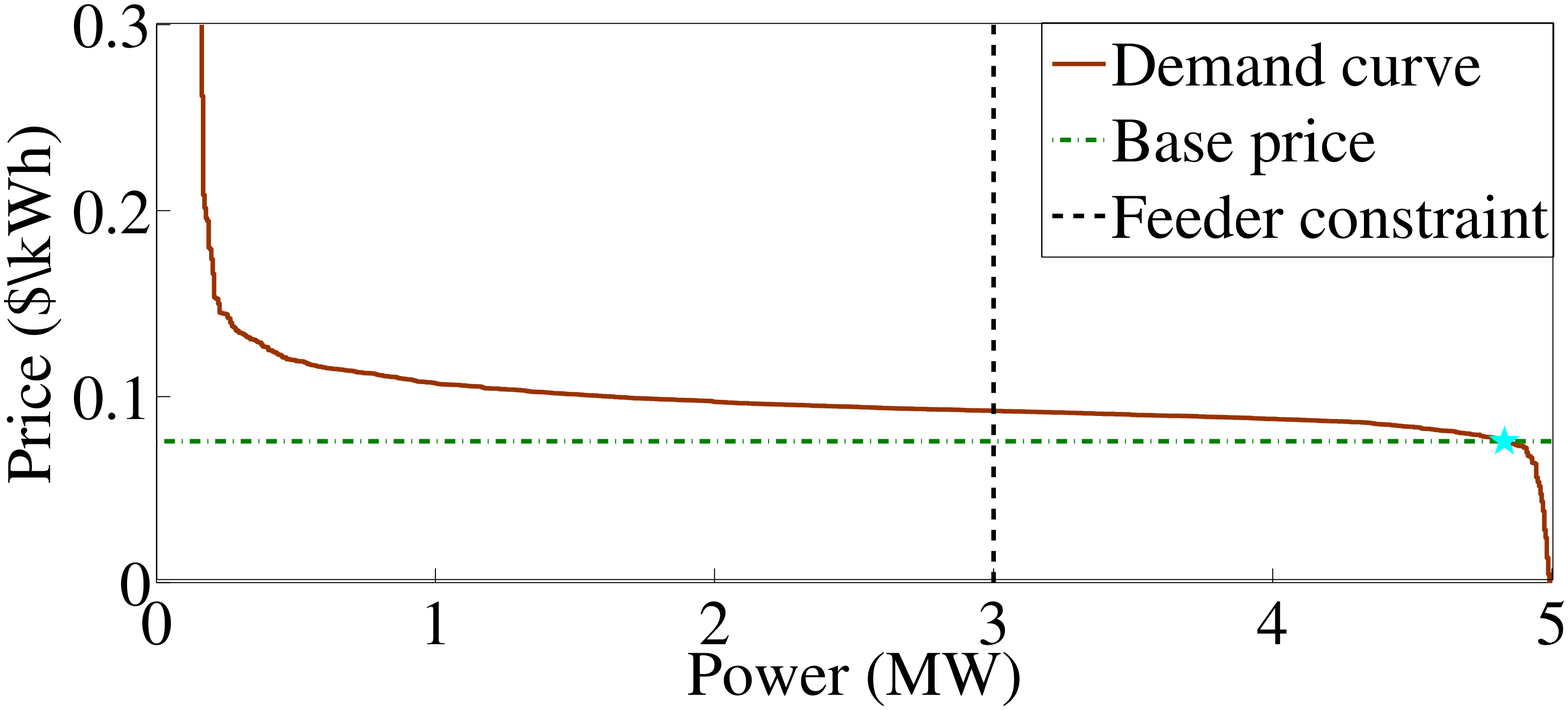}
\caption{The market clearing strategy during power congestion under real time pricing. The market clearing point violates the feeder capacity constraint.}
\label{con_rtp}
\end{figure}

\begin{figure}[t]
\centering
\includegraphics[width = 0.85\linewidth]{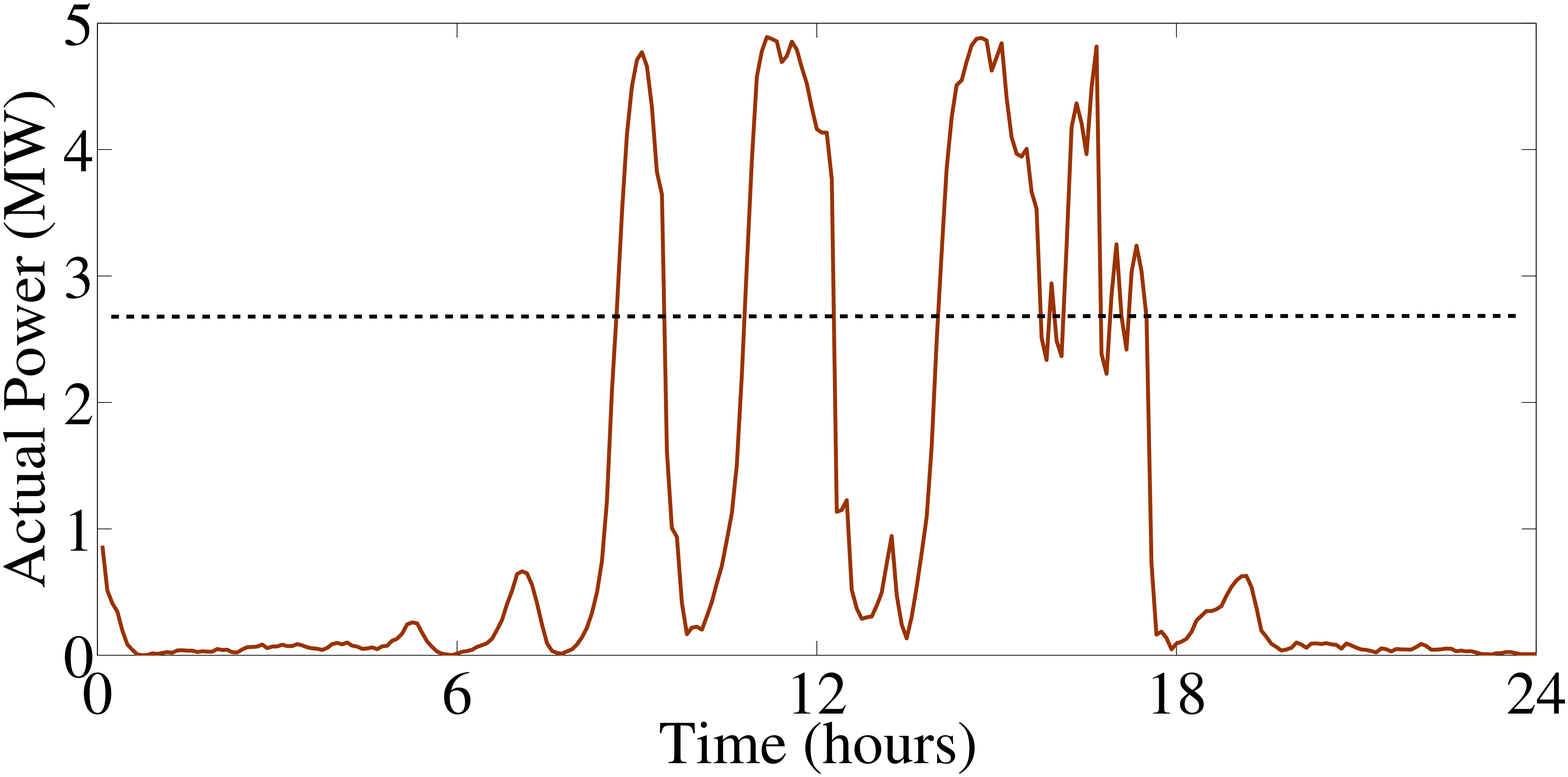}
\caption{The actual power trajectory and the feeder capacity. The outside air temperature record is on August 16, 2009 in Columbus, OH.}
\label{gridwise}
\end{figure}

\begin{figure}[t]
\centering
\includegraphics[width = 0.85\linewidth]{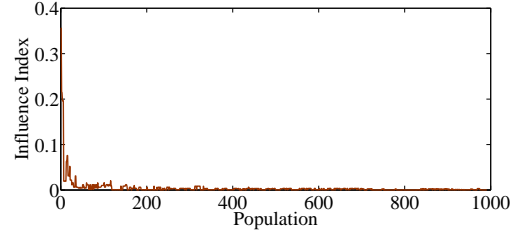}
\caption{The influence index of a randomly chosen user over different population sizes. The influence of individual bid drops rapidly as the population size grows.}
\label{influence}
\end{figure}

\begin{figure}[t]
\centering
\includegraphics[width = 0.80\linewidth]{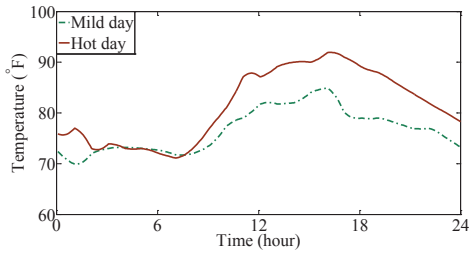}
\caption{The hourly temperature used in the simulation. The record is from August 16 (hot day) and August 20 (mild day), 2009 in Columbus, OH.}
\label{outsideair}
\end{figure}

\subsection{Comparison with Other Strategies}
In this subsection we first compare the proposed mechanism with RTP \cite{allcott2009real}. RTP can incentivize users to shift demand from high price periods to low price periods to reduce electricity expenditures. However, as such an approach directly passes the base price to the retail market, it can not achieve predictable and reliable aggregated power response, which is essential in many demand response programs. To illustrate these limitations, we compare our framework with RTP by applying RTP in the considered problem. In this simulation the coordinator clears the market by directly passing the base price to individual users, and the devices respond to the energy price according to the response curve described in Fig. \ref{fig:response}. Except for the pricing strategy, all the parameters are the same as in the simulation in Section III-B, and the result is presented in Fig. \ref{RTP}.
When there is no congestion, the real time pricing scheme has the same performance as the proposed mechanism, and efficient energy allocation can be achieved. However, during the power congestion period, the RTP method can not prevent the aggregated power from exceeding the feeder power limit. For instance, the market clearing process at 4:40 PM is presented in 
Fig. \ref{con_rtp}. In the example, due to the increased power demand in the afternoon, the market clearing point exceeds the feeder power limit. This issue can be solved in the proposed mechanism with an elevated energy price during power congestion. We emphasize that the proposed mechanism may also fail to cap the energy in certain extreme case. For instance, when the outside temperature is extremely high and all the participating TCLs have very small thermal capacity and resistance, then it is possible that a large number of users have to turn on the TCL for the entire market period, and the  aggregated energy can not be capped effectively. However, this is not due to the proposed mechanism, but rather because of the physical limitation of the system.

In addition, the proposed mechanism is also compared with the original Gridwise${\textsuperscript{\textregistered}}$ demonstration project (base scenario). In the simulation, the market clearing strategies of the two cases are the same, while the bidding strategies are different. In the base scenario, each device submits a bid based on the current temperature, while the device in the proposed mechanism computes the bid according to Fig. \ref{biddingana}. Except for the bidding strategy, all the parameters are the same as in the simulation in Section III-B, and the result of the base scenario is presented in Fig. \ref{gridwise}. In this case, the user bids only depend on the room temperature, and the information regarding the model parameter is missing. Therefore, although the same pricing strategy is applied, the coordinator still can not achieve the desired aggregated power response.

\subsection{Impact of Some Key Parameters}
This subsection discusses how a few important parameters can affect the performance of the proposed mechanism. 
\subsubsection{Number of Households}
In this paper, we assume that every user is a price taker: the bid of an individual user can not affect the market price. Theoretically, this assumption only holds when the market satisfies some regulatory conditions, such as sufficiently many users, free entry, homogeneous good, etc \cite[Chap. 12]{mas1995microeconomic}. In this subsection we use numerical simulations to investigate to what extend this assumption can be justified. In particular, we simulate the influence of individual bid on the market price, and explore how this influence changes with the growing number of participating households.  This can be done by perturbing the bidding price of a user $i$, and see how the market price changes with this perturbation while the bids of all the other users remain the same. It can be verified that under the proposed clearing strategy, the user bid could only affect the market price in two possible ways: change the market price by a fixed value (regardless of how big the perturbation is) or no influence at all. To quantitatively represent this market influence, we define an \emph{influence index}, which is the maximum market price change (in percentage) that a user could incur by perturbing its bidding price. When the individual bid has no influence on the market clearing price, the influence index is zero, and the price-taker assumption holds.

\begin{figure}[t]
\centering
\includegraphics[width = 0.85\linewidth]{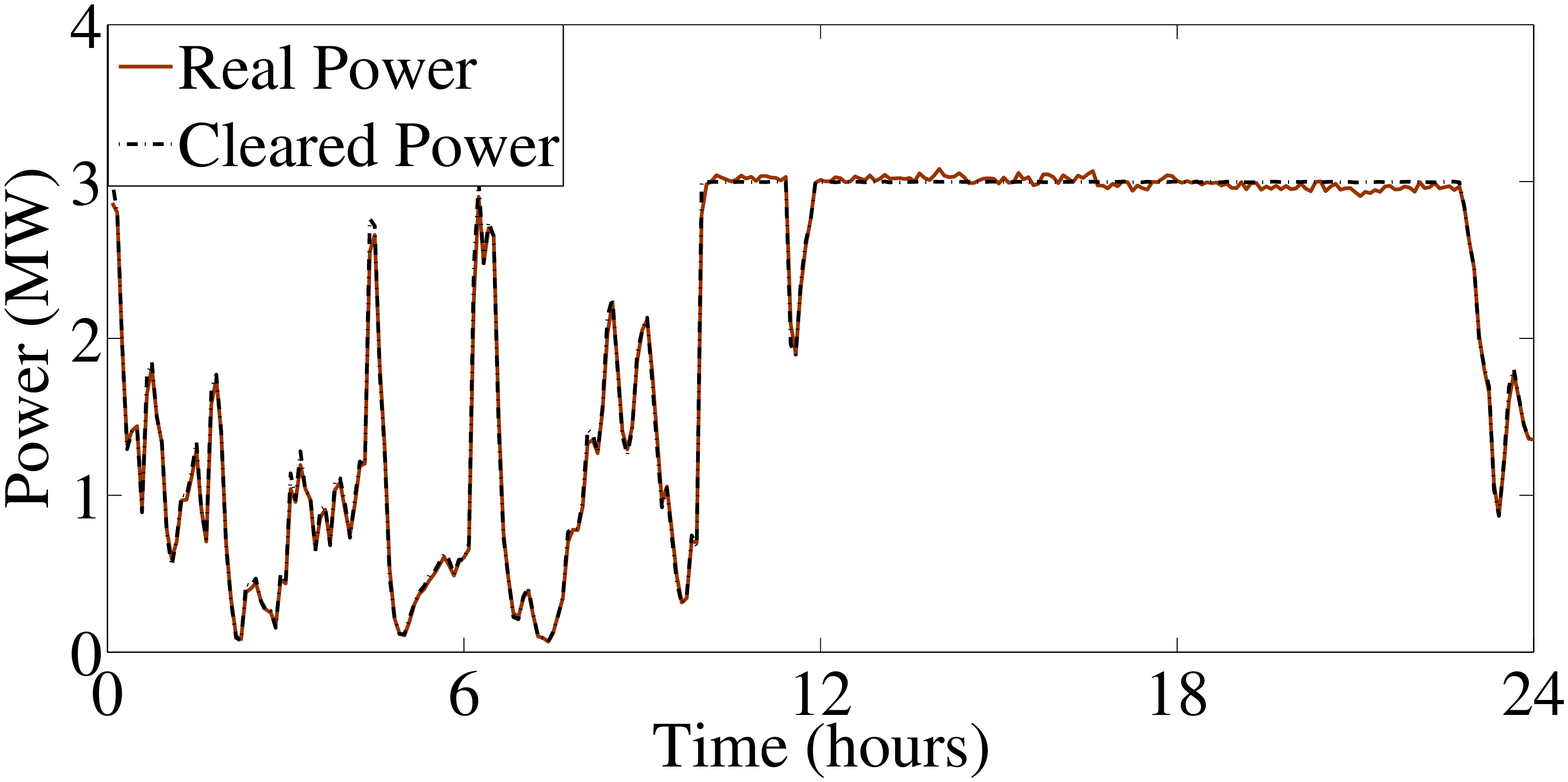}
\caption{Comparison of the actual power trajectory and the cleared power. The outside air temperature record is on August 16, 2009 in Columbus, OH.}
\label{powertrajectory2}
\end{figure}

\begin{figure}[t]
\centering
\includegraphics[width = 0.85\linewidth]{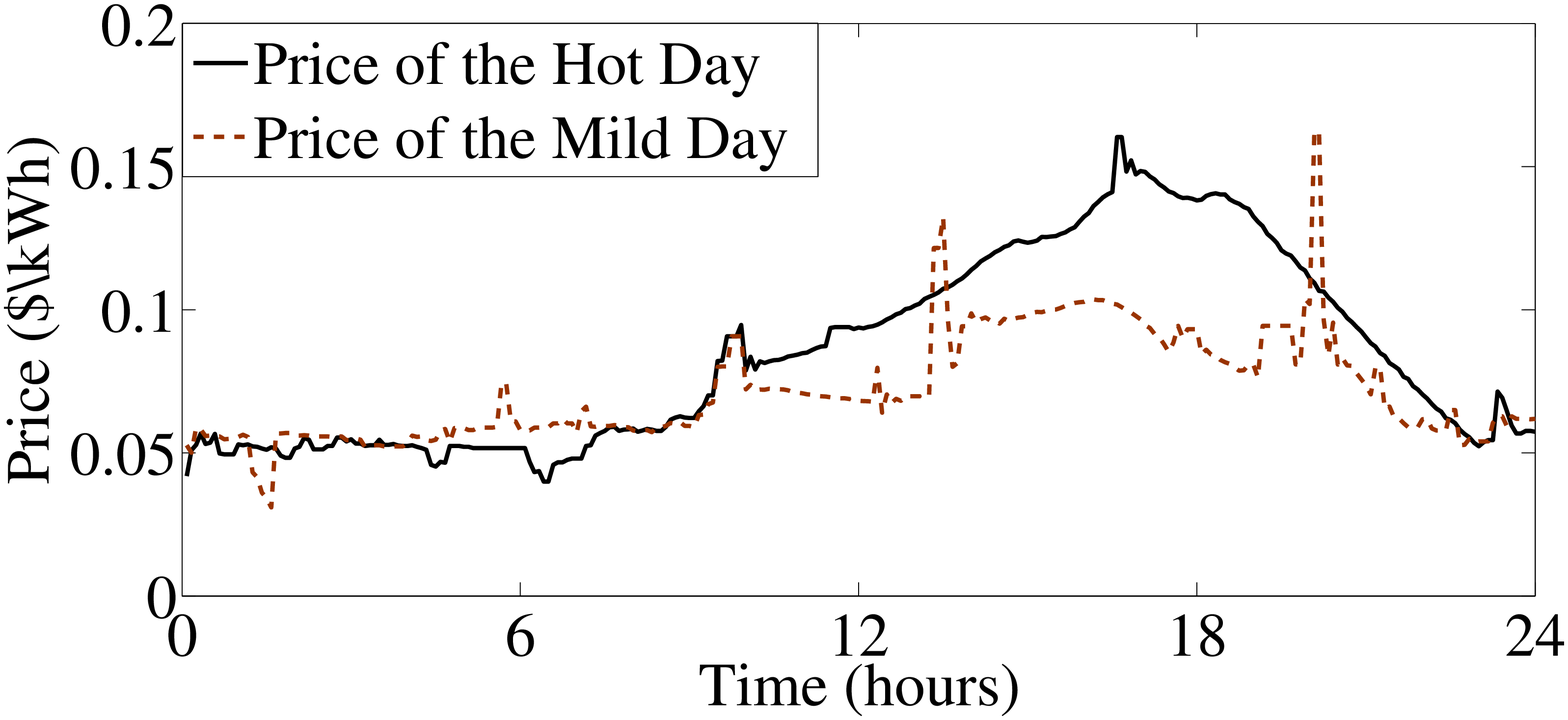}
\caption{The market prices of the hot day and the mild day. The average price of the hot day is higher than that of the mild day.}
\label{price_hot}
\end{figure}

\begin{figure}[t]
\centering
\includegraphics[width = 0.8\linewidth]{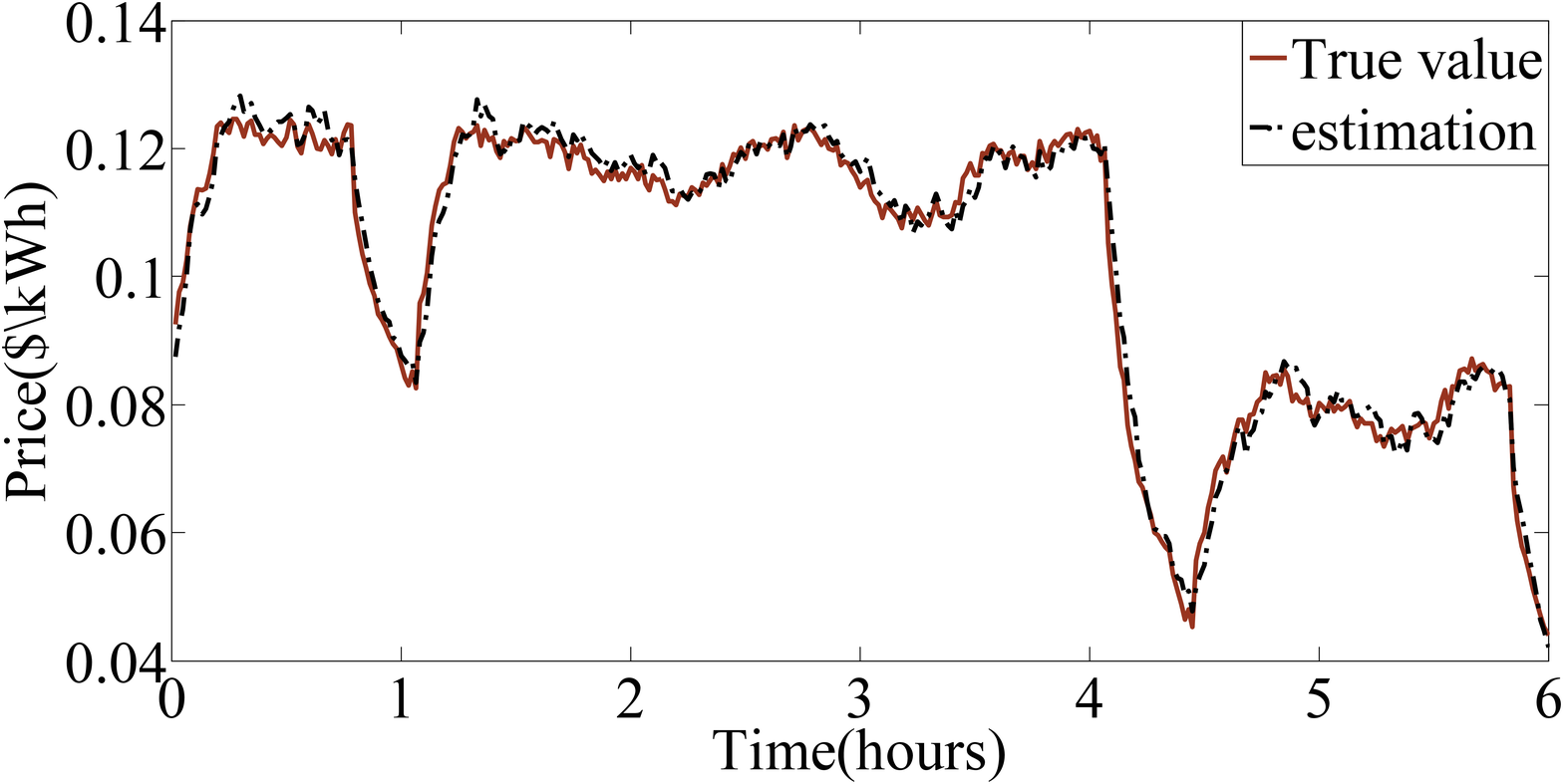}
\caption{ The estimation result of the output-based bidding algorithm when the initial guess is randomly selected from $50\%$ and $150\%$ of its true value.}
\label{EM2}
\end{figure}

\begin{figure}[t]
\centering
\includegraphics[width = 0.75\linewidth]{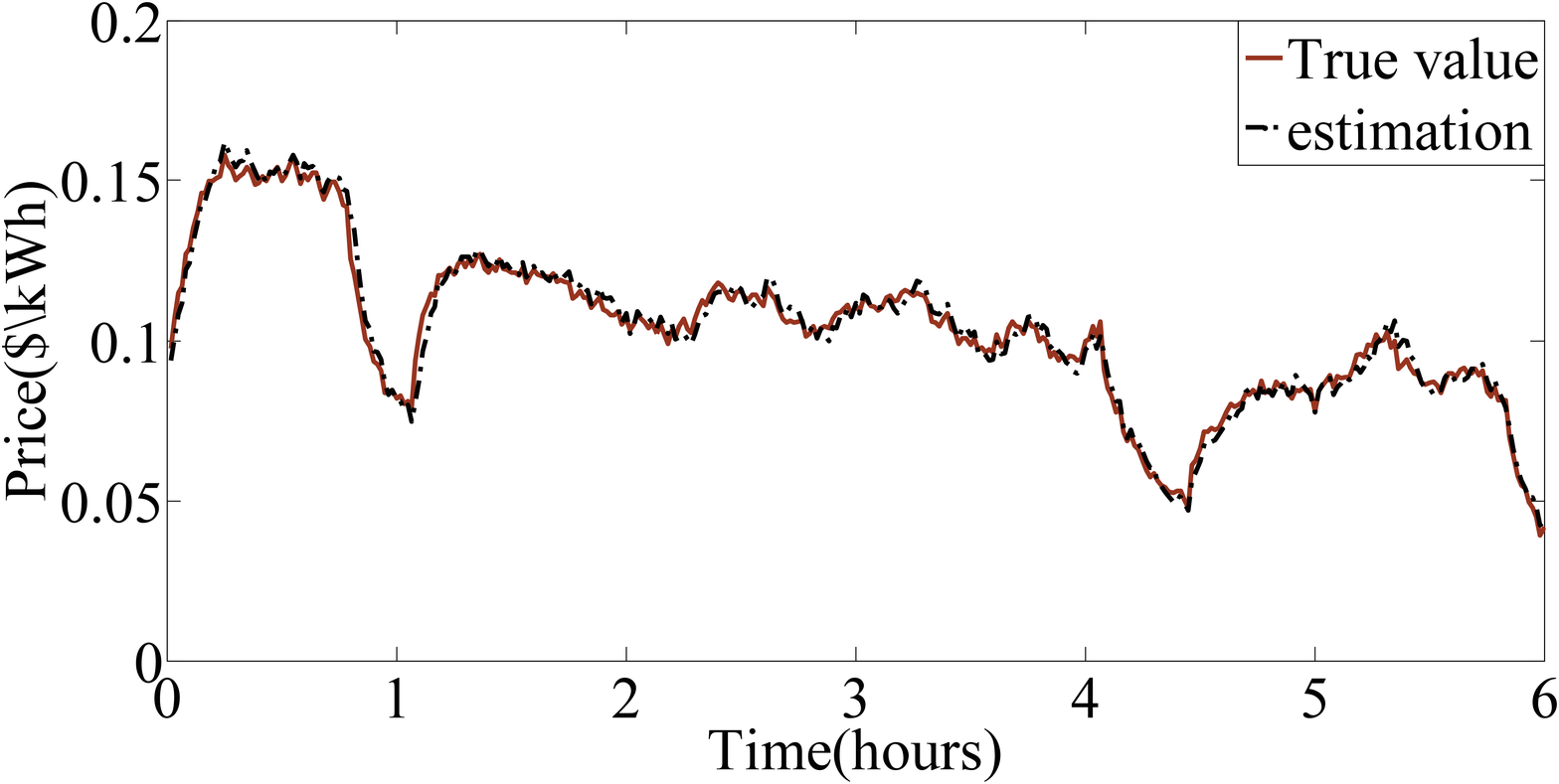}
\caption{The estimation result of the output-based bidding algorithm when the process noise and measurement noise are generated based on a uniform distribution.}
\label{EM3}
\end{figure}

The simulation can be done in the following steps. First, randomly choose a group of users (for example, 100 users) and a market period, simulate the market bidding and clearing process, and derive the corresponding market clearing price. Second, choose one user from this group, perturb his bid, and rerun the market clearing process to obtain another market price. Third, compute the influence index based on the market clearing prices derived from the first two steps. Fourth, enlarge the group and repeat all the procedures described above. Notice that when there is no congestion (as illustrated in Fig. \ref{clearnocon}), the clearing price is always the base price, and the influence index is 0. Therefore, all the simulations are done in a market period during which power congestion occurs. In addition, to enforce a fair comparison, we assume that the feeder capacity constraint changes according to the number of participating household. For example, if the maximum power of each air conditioning load is 5kW, and there are $N$ loads in the project, then the maximum aggregated power is $5N$ kW, and the feeder power capacity is $60\%$ of the maximum power, i.e., $3N$ kW. 

The simulation result is shown in Fig. \ref{influence}. It starts with 10 users and the inference index in this case is around $35\%$. This influence drops rapidly as the number of the participating loads increases. When there are more than 200 loads, the influence index is always less than $1\%$. When the population size is larger than 500, the inference index is less than $0.4\%$, in which case the influence of individual users on the market price can be safely neglected.

\subsubsection{Weather Information}
Aside from the number of participating households, the outside temperature data is also an important parameter that affects the performance of the proposed mechanism. The high temperature period can significantly increase the aggregated power demand of the air conditioning loads, and therefore cause more power congestion. For this reason, we evaluate the proposed method with a different temperature record. The data is obtained from \cite{weatherdata} on August 16 (hot day), Columbus, OH, as shown in Fig. \ref{outsideair}.

The power trajectory and the market clearing prices are presented in Fig. \ref{powertrajectory2} and Fig. \ref{price_hot}, respectively. Since the elevated temperature increases the power demand, much more power congestion can be observed during the hot day. Despite the power congestion, the simulation result shows that the proposed framework can still effectively cap the aggregated power below the power limit, and the actual power accurately matches the planned power.

\subsubsection{Initial Guess of the output-based bidding algorithm}
The initial guess of the output-based bidding algorithm is also crucial to performance of the estimation result. In our previous simulations, the initial guess is generated by randomly selecting a value between $90\%$ and $110\%$ of its true value. Therefore, to implement the output-based bidding algorithm, we need to assume that users have some prior knowledge of the unknown parameters to guarantee that initial guess is within this range (from $90\%$ to $110\%$ of the true value). In this subsection we explore to what extend we can relax this assumption without compromising the estimation performance. In particular, we use the same model parameters as in Section III-C, and test the proposed algorithm with an error of $50\%$. The estimation result is shown in Fig. \ref{EM2}, which shows that the output-based bidding algorithm can accurately estimate the bidding prices even with $50\%$ error on the initial guess.

\subsubsection{Model noises}
The proposed EM algorithm is developed mainly based on the linear Gaussian model (\ref{eq:uncertainmodel}), where we assume that $B_{on}^i$ and $B_{off}^i$ can be decomposed into three parts: the external signals, a constant and a Gaussian process noise. However, we emphasize that in our particular problem, the process noise does not have to be Gaussian, and the proposed EM algorithm can be extended to a much broader class of ``real'' dynamical systems. To support this argument, the proposed EM algorithm is tested with non-Gaussian noises. In particular, in the simulation the process noises and the measurement noises are generated based on a uniform distribution, while the rest of the model parameters are generated the same way as described in Section III-C. The simulation result is presented in Fig. \ref{EM3}, where the estimated bidding prices are close to the real bidding prices. 
The key reason for this result is that in our problem, the EM algorithm does not need to accurately estimate all the unknown parameters (such as $\bar{A}_i, \bar{B}_i$ and $\bar{C}_i$). Instead, we only need to estimate the bidding price, which is a scalar-valued function of these unknown parameters. According to  (\ref{eq:bidding}), the bidding price mainly depends on the current room temperature $T_c^i(t_k)$ and the room temperature 5 minutes after $t_k$. Since we have measurements of the room temperature for the past 6 hours in the algorithm, this information is already contained in these measurements except for the last 5 minutes of the 6 hour period. Therefore, although the algorithm can not converge to all the true unknown parameters, it does converge to the true bidding price under non-Gaussian noise distributions.

\section{conclusion}
This paper presents a market mechanism for the coordination of thermostatically controlled loads, where a coordinator manages a group of TCLs using pricing incentives to maximize the social welfare subject to a peak energy constraint.  In the paper, a mechanism is proposed to implement the desired social choice function in dominant strategy equilibrium. This mechanism consists of a novel bidding strategy that incorporates information on both the load dynamics and the time-varying user preferences. It is proven that under the proposed mechanism, the coordinator can not only maximize the social welfare but also realize the team optimal solution. Future work includes formulating the fully dynamic market-based coordination framework with multiple periods and extending the results to energy storage devices and deferrable loads such as plug-in electric vehicles, washers, dryers, among others.

\appendix{}
\subsection{Proof of Proposition 1}
When each device submits $h_i$ as the bid, we have $b_i(\cdot;\theta_i)=h_i(\cdot;\theta_i)$. According to (\ref{EQ_1}), each user will receive an energy allocation that satisfies $a_i^*=h_{i}(P_c;\theta_i)$. Based on (\ref{eq:userresponse}), we have: $a_i^*=\argmax_{0\leq a_{i} \leq E_i^{m}} V_i(a_{i};\theta_i)-P_ca_{i}$. Therefore, when $b_i(\cdot;\theta_i)=h_i(\cdot;\theta_i)$, the resulting energy allocation maximizes the utility of each user. According to Definition 2, the strategy profile $(h_1(\cdot;\theta_1),\ldots,h_N(\cdot;\theta_N))$ is a dominant strategy equilibrium of the proposed mechanism. 

\subsection{Proof of Proposition 2}
Notice that the social choice function characterizes the optimal solution to the coordinator's optimization problem (\ref{eq:coordinatoropt}), and the team solution provides an upper bound on the social welfare for (\ref{eq:coordinatoropt}). Therefore, to prove Proposition 2, it is sufficient to show that the proposed pricing strategy realizes the team solution. 

Based on Proposition 1, $b_i=h_i$. Therefore, we have the following relations:
\begin{align}
\begin{cases}
a_i^*=h_i(P_c^*;\theta_i), \text{ for all } i=1,\ldots,N \\
P_c^*=\text{max} \big \{\bar{P},P^* \big \}\\
P^*=C'\left( \sum_{i=1}\nolimits ^N a_i^*\right) \\
\sum_i^N h_i(\bar{P}, \theta)=D.
\end{cases}
\label{eq:optpricing}
\end{align}
In addition, the KKT condition for the $i$th user's individual utility maximization problem (\ref{eq:individualopt}) is as follows:
\begin{equation}
\label{kkt}
-V'_i(a_i^*;\theta_i)+P_c^*+u_1^i-u_2^i=0,
\end{equation}
where $u_1^i$ and $u_2^i$ are the Lagrangian multiplier satisfying:
\begin{align}
\begin{cases}
u_1^i\geq 0, u_2^i\geq 0 \\
u_1^i=0 \quad \text{if} \quad a_i^*\neq E_i^m\\
u_2^i=0 \quad \text{if} \quad a_i^* \neq 0.
\end{cases}
\label{eq:lagrangian}
\end{align}
Define $u=P_c^*-C'\left( \sum_{i=1}\nolimits ^N a_i^*\right)$, then equation (\ref{kkt}) becomes:
\begin{equation}
\label{kkt1}
-V'_i(a_i^*;\theta_i)+C'\left( \sum_{i=1}\nolimits ^N a_i^*\right)+u+u_1^i-u_2^i=0,
\end{equation}
According to (\ref{eq:optpricing}), when $\sum_i^N a_i^*< D$, we have $P_c^*=P^*=C'\left( \sum_{i=1}\nolimits ^N a_i^*\right)$, therefore, $u=0$. When $\sum_i^N a_i^*=D$, we have $P_c^*=\bar{P}$, and therefore, $u=\bar{p}-p^*$. Since $h_i$ is non-increasing, we have $u\geq 0$. This indicates that $u$, $u_1^i$ and $u_2^i$ are the Lagrangian multipliers of the team problem, and (\ref{kkt1}) is exactly the KKT condition for the team problem (\ref{eq:teamopt}). Since the team problem is a concave optimization problem, the KKT conditions are also sufficient. Thus $a^*=(a_1^*,\ldots,a_N^*)$ is the team solution. This completes the proof. 

{
\bibliographystyle{unsrt}
\bibliography{TransactiveControl}
}

\end{document}